\documentclass[12pt,a4paper]{article}
\usepackage{amsmath,amssymb}
\usepackage{latexsym}
\usepackage{graphicx}
\usepackage{color}
\usepackage{multirow}
\usepackage{indentfirst}
\usepackage{subfigure}
\usepackage{mathrsfs}
\usepackage[pdftex]{hyperref}
 \usepackage[numbers,sort&compress]{natbib}
\usepackage{epstopdf}
\usepackage{epsfig}

\newtheorem{exam}{\hspace{6mm}Example}[section]

\begin{document}
\baselineskip=2pc

\begin{center}
{\large \bf  A compact simple HWENO scheme with ADER time discretization for hyperbolic conservation laws II: triangular meshes}
\end{center}

\centerline{
Dongmi Luo%
\footnote{Institute of Applied Physics and Computational Mathematics, Beijing 100088, China. E-mail: dongmiluo@stu.xmu.edu.cn.},
Zhuang Zhao*%
\footnote{School of Mathematical Sciences and Fujian Provincial
	Key Laboratory of Mathematical Modeling and High-Performance
	Scientific Computing, Xiamen University, Xiamen, Fujian 361005, China. E-mail: zzhao@xmu.edu.cn.},
Jianxian Qiu%
\footnote{School of Mathematical Sciences and Fujian Provincial
Key Laboratory of Mathematical Modeling and High-Performance
Scientific Computing, Xiamen University, Xiamen, Fujian 361005, China. E-mail: jxqiu@xmu.edu.cn.},
Yibing Chen%
\footnote{Institute of Applied Physics and Computational Mathematics, Beijing 100088, China, 
E-mail: chen\_yibing@iapcm.ac.cn.}

}
\vspace{20pt}

\begin{abstract}
A compact and high order HWENO scheme using ADER (Arbitrary high order using DERivatives) time discretization is developed for hyperbolic conservation laws on the triangular mesh, which is the extension of the work on the structured mesh (Luo et. al. (2024) \cite{luo2023}). The Lax-Wendroff procedure is employed to convert time derivatives to spatial derivatives. Thanks to this, the cell averages of the derivatives of the solution can be obtained by the time accurate solution as Gaussian points along the cell interfaces through the Green-Gauss theorem instead of by the evolution solution directly in the conventional HWENO methods. Comparing with the existing Runge-Kutta HWENO (RK-HWENO) method on the unstructured mesh (Zhao et. al. (2025) \cite{zhao2025}), the new method has the following advantages. Firstly, the RK-HWENO method must solve the additional equations for reconstructions and time advancing, which is avoided for the new method. Secondly, the HWENO reconstruction in the new method is performed once per time step and is different from the RK-HWENO method, in which the reconstruction is performed several times every time step. Because of these advantages the new method is more efficient than the RK-HWENO method with smaller numerical errors and less computational costs. Besides, comparing with the existing ADER-WENO methods \cite{dumbser20071,dumbser20072} under the same order of accuracy, the stencil of the new method is more compact since the both the function and its first derivative values are used in the reconstruction of the HWENO schemes. Numerical examples demonstrate that the new method can achieve the high order for smooth solutions both in space and time, keep non-oscillatory near discontinuities.
\end{abstract}


\pagenumbering{arabic}

\newpage

\section{Introduction}
\label{sec1}
\setcounter{equation}{0}
\setcounter{figure}{0}
\setcounter{table}{0}

In \cite{luo2023}, a compact simple Hermite  weighted essentially non-oscillatory (HWENO) scheme with ADER (Arbitrary high order using DERivatives)  time discretization for hyperbolic conservation laws was proposed on the structured meshes. 
In this paper, we go on with the work to develop the compact simple HWENO scheme for the numerical solutions of hyperbolic conservation laws in the form 
\begin{equation}
\label{eq2}
\left\{
\begin{array}{l}
W_t + \nabla\cdot F=0,\\
W(\textbf{x},0)=W_0(\textbf{x}),
\end{array}
\right.
\end{equation}
on the triangular meshes,
where $\textbf{x}=(x,y)\in \Omega \subset \mathbb{R}^2$, $\Omega$ is a bounded domain, $F=(f(W),g(W))$, and $W, f(W), g(W)$ are either scalars or vectors. 

The HWENO schemes are first proposed for solving one-dimensional nonlinear hyperbolic conservation law systems \cite{qiu2003}, whose idea of the reconstruction comes from the original WENO schemes \cite{liu1994,jiang1996}. However both the function and its derivative values are evolved in time and employed in the reconstructions while the function values are adopted in the original WENO schemes \cite{liu1994,jiang1996}. Then the HWENO schemes are extended to the two-dimensional structured meshes \cite{qiu2005} and unstructured meshes \cite{zhu2009}. But these methods have some drawbacks. Firstly, the linear weights of the point values in two dimensions are obtained by the least square method. Secondly, the reconstructions are performed twice with two different stencils. In addition, the additional equations need to be solved for time evolution.

To deal with the first problem, a new type of WENO methods for hyperbolic conservation laws were developed \cite{zhu2016,zhu2017,zhu2018}, which has the advantages that the associated linear weights can be any positive numbers with the only requirement that their summation equals one. Then the idea is extended to the multi-resolution WENO \cite{zhu2018m,zhu2019} and HWENO schemes \cite{lij2021,li2022,li2021,fan2023,zhao2025}. To avoid the second problem, a fifth-order HWENO scheme with unified stencils was proposed for hyperbolic conservation laws \cite{fan2025}. The key idea of the scheme is to modify the first-order moment by a HWENO limiter only in the time discretizations using the same spatial reconstructions. Thus, the HWENO scheme needs only one HWENO reconstruction with one set of stencils. However, the additional equations are required to be solved for time evolutions in all the HWENO schemes mentioned above. To address the problems, a fifth order accurate HWENO reconstruction method for hyperbolic conservation laws in the framework of two-stage fourth order accurate temporal discretization was developed \cite{du2018}. They used the interface values, which are available using the generalized Riemann problem (GRP) solver, and avoided computing the derivatives of the solution for time evolution in the conventional HWENO  \cite{qiu2003,qiu2005,zhu2009}. However, the reconstruction needs to be performed several times for the multi-stage time discretization methods. To avoid the disadvantages of the multi-stage temporal discretization methods, many high order one-stage method such GRP, HGKS (High order Gas-Kinetic Scheme) \cite{li2021,ji2020,ji2018}, ADER and other Lax-Wendroff procedure based methods \cite{qiu2003lw,qiu2005lw,qiu2007lw} are proposed.

The GKS has been systematically developed in recent years \cite{ji2018,ji2020,li2021}. A fourth-order and third-order HWENO reconstruction based on the HGKS is proposed on the structured mesh \cite{ji2018} and the unstructured mesh \cite{ji2020}, respectively. Recently, a compact and efficient HGKS (CEHGKS) \cite{li2021} based on a HWENO reconstruction using the idea is proposed for Euler equations on the structured mesh, which is based on the framework of a one-stage efficient HGKS and can achieve high order accuracy both in space and time. And the method overcomes the drawbacks of the classic HWENO scheme. However, the technique is only suitable in the framework of GKS.

 Another popular one-stage temporal discretization method is ADER method, which was first put forward by Toro and collaborators for linear problems on Cartesian meshes \cite{toro2001}. Soon after that the ADER methods are extended to many nonlinear problems \cite{titarev2002, titarev2005, toro2006,dumbser20071,dumbser20072,kaser2005,boscheri2016,boscheri2016caf,castro2012}. However, most of these ADER method are combined with the classic WENO methodology and the stencil used in the reconstructions is becoming wider with an increasing the order of accuracy, which makes the method more complex on high dimension and unstructured meshes \cite{boscheri2016,boscheri2016caf,kaser2005,castro2012}.  To avoid the disadvantages of ADER-WENO methods, the ADER approaches are extended to nonlinear systems in the framework of DG methods (ADER-DG)\cite{dumbser2006,fambri2017}. However, ADER-DG methods also need to solve the additional equations for obtaining the internal freedom for reconstructions and time advancing, which increases the computational costs. Recently, a compact and simple HWENO method using the ADER time discretization is proposed for hyperbolic conservation laws on the structured meshes \cite{luo2023}. 

 To accommodate computations involving the complex geometries,
 a third order HWENO method for hyperbolic conservation laws based on the ADER approach is proposed on the triangular mesh in this paper, which is denoted by ADER-HWENO and an extension of the work \cite{luo2023}. Different from the conventional HWENO methods \cite{qiu2003,qiu2005,zhu2009}, the cell averages of the derivatives of the solution are updated by the time accurate solutions as Gaussian points along the cell interfaces through the Green-Gauss theorem instead of by the evolution solution directly, which are also used in \cite{ji2020,li2021,luo2023}. Owing to this, the HWENO reconstruction is performed once with one stencil and avoids the additional equation for time evolution the same as in \cite{luo2023}. Thus the new method is more efficient with smaller numerical errors and less computational costs than the HWENO scheme in \cite{zhao2025}, in which the Runge-Kutta method is employed for the temporal discretization and the additional equations need to be solved for the HWENO reconstruction and time advancing. In addition, comparing with the existing ADER-WENO methods \cite{dumbser20071,dumbser20072} at the same order of accuracy, the stencil of the new method is more compact since both the function and its derivative values are adopted in the reconstruction of the HWENO schemes.  Moreover, the new method makes the best use of the information in the ADER method, which leads that the time evolution of the cell averages of the derivatives is simple. Numerical tests show that the new method can achieve high order for smooth solutions in space and time and keep non-oscillatory near discontinuities.
  
 An outline of the paper is given as follows. The HWENO scheme based on the ADER method is described in Section \ref{sectwo} for the hyperbolic conservation laws  on the triangular mesh in detail.
 A selection of numerical examples is presented to demonstrate the accuracy and the capability of the ADER-HWENO method in Section \ref{secnumexam}. In Section \ref{secconclusion}, conclusions are drawn. 
 
\section{The numerical scheme on the triangular meshes}
\label{sectwo}
\setcounter{equation}{0}
\setcounter{figure}{0}
\setcounter{table}{0}
In this section, we describe the ADER-HWENO method for two dimensional problems on the unstructured meshes. We consider the hyperbolic conservation laws \eqref{eq2}.
We assume a triangular mesh $\mathscr{T}_h$ for the domain $\Omega$ is given and $K\in \mathscr{T}_h$. 

Integrating \eqref{eq2} on $K\times (t_n,t_{n+1})$, we obtain
\begin{align*}
\Bar W(\text{\bf {x}}_K,t_{n+1})=\Bar W(\text{\bf {x}}_K,t_{n})-&\frac{\Delta t}{|K|}(\frac{1}{\Delta t}\int_{t_n}^{t_{n+1}}\int_{\partial K}F\cdot \vec{n}dsdt),
\end{align*}
where $\Bar W(\text{\bf {x}}_K,t_{n})=\frac{1}{|K|}\int_{K}W(\text{\bf {x}},t_n)dxdy$, $\Delta t=t_{n+1}-t_n$, $\vec{n}=(n_x,n_y)^T$ is the outward unit normal vector of the triangular boundary $\partial K$ and $|K|$ is the volume of the element $K$. Then the finite volume method in two dimensions is given by
\begin{equation}
\label{flx}
\Bar W_{K}^{n+1}=\Bar W_{K}^{n}-\frac{\Delta t}{|K|}\Hat F_{K},
\end{equation}
where $\Bar W^n_K$ is a high order approximation to $\Bar W(\text{\bf {x}}_K,t_{n})$, and $\Hat F_{K}\approx \frac{1}{\Delta t}\int_{t_n}^{t_{n+1}}\int_{\partial K}F(W(\text{\bf{x}},t))\cdot\vec{n}dsdt.$ 
 One can observe that this leads to equations including the line integrals on the cell $K$, which can be computed by the Gaussian quadrature rule
\begin{align}
\label{fx}
\Hat F_{K}&=\sum\limits_e\sum\limits_{G_e}(\frac{1}{\Delta t}\int_{t_n}^{t_{n+1}}F(W(\text{\bf{x}}_{G_e},t))\cdot\vec{n}dt)w_{G_e}|e|,
\end{align}
where $e$ represents the edges of the element $K$, $\text{\bf{x}}_{G_e}$ and $w_{G_e}$ represent the Gaussian points on $e$ and the corresponding weights, respectively. The summations $\sum\limits_e$ and $\sum\limits_{G_e}$ are taken over the edges of $\partial K$ and Gauss points on $e$, respectively. In order to obtain the high order of accuracy in both time and space, the ADER approach \cite{luo2023,titarev2002} is employed. 

\subsection{ADER method on unstructured meshes}
For simplicity, the idea of the paper \cite{luo2023} for Navier-Stokes equations is adopted. Therefore, to evaluate the numerical flux an appropriate Gaussian rule for time integration is used in \eqref{fx}
\begin{align}
	\label{flux}
	 \hat F_{K}&=\sum\limits_e\sum\limits_{G_e}\sum\limits_{s=1}^{K_t}(F(W(\text{\bf{x}}_{G_e},\lambda_s\Delta t))\cdot\vec{n}w_{s})w_{G_e}|e|=\sum\limits_e\sum\limits_{G_e}\sum\limits_{s=1}^{K_t}\hat Hw_{s}w_{G_e}|e|,
\end{align}
where $\lambda_s$ and $w_{s}$ are properly scaled nodes and weights of the rule, $K_t$ is the number of the nodes, and $\Hat H$ is the numerical flux. In this paper, the simple local Lax-Friedrichs flux is employed for the advection term, i.e.,
\[
\hat H=\frac{1}{2}(H(W^{ext}(\text{\bf x}_{G_e},\tau)) + H(W^{int}(\text{\bf x}_{G_e},\tau))-\alpha (W^{ext}(\text{\bf x}_{G_e},\tau)-W^{int}(\text{\bf x}_{G_e},\tau))),
\]
where $\tau=\lambda_\alpha\Delta t, H(W^{ext}(\text{\bf x}_{G_e},\tau))=F(W^{ext}(\text{\bf x}_{G_e},\tau))\cdot\vec{n}, H(W^{int}(\text{\bf x}_{G_e},\tau))=F(W^{int}(\text{\bf x}_{G_e},\tau))\cdot\vec{n}$, 
and $W^{int}(\text{\bf{x}}_{G_e})$  and $W^{ext}(\text{\bf{x}}_{G_e})$ are defined as the values from the interior and exterior of $K$, respectively, i.e.,
\[
W^{int}(\text{\bf{x}}_{G_e},t)=\lim\limits_{\text{\bf{x}}\to\text{\bf{x}}_{G_e},\text{\bf{x}}\in K}W(\text{\bf{x}},t),
\quad W^{ext}(\text{\bf{x}}_{G_e},t)=\lim\limits_{\text{\bf{x}}\to\text{\bf{x}}_{G_e},\text{\bf{x}}\notin K}W(\text{\bf{x}},t),
\]
and $\alpha$ is the numerical viscosity constant taken as the largest eigenvalues in magnitude of
\[
\frac{\partial}{\partial W}(F(\bar W_K)\cdot\vec{n}),\quad \frac{\partial}{\partial W}(F(\bar W_{K'})\cdot\vec{n}),
\]
where $K$ and $K'$ are the elements sharing the common edge $e$.

Then the approximate solutions at time $t=\tau$ on each side of the interface are evaluated by
\begin{align}
	\label{tay}
W^{int}(\textbf{x}_{G_e},\tau)&\approx W^{int}(\textbf{x}_{G_e},0^+) + \sum\limits_{k=1}^2[\partial_t^kW^{int}(\textbf{x}_{G_e},0^+)]\frac{\tau^k}{k!},\\
\label{tayext}
W^{ext}(\textbf{x}_{G_e},\tau)&\approx W^{ext}(\textbf{x}_{G_e},0^+) + \sum\limits_{k=1}^2\partial_t^kW^{ext}(\textbf{x}_{G_e},0^+)]\frac{\tau^k}{k!},
\end{align}
where
\[
\partial^{(k)}_tW(\text{\bf{x}},t)=\frac{\partial^k}{\partial t^k}W(\text{\bf{x}},t),\quad 0^+\equiv \lim_{\tau\rightarrow 0^+}\tau.
\]
Similar to the classic ADER method, all the time derivatives are also expressed as functions of space derivatives by the Lax-Wendroff procedure based the original equations. Then one can get the following expressions:
\begin{align*}
\partial_t W &= -(\frac{\partial f}{\partial W})\partial_x W - (\frac{\partial g}{\partial W})\partial_y W,\\
\partial_{tx}W &= -(\frac{\partial^2 f}{\partial W^2}\partial_x W)(\partial_x W)  -(\frac{\partial f}{\partial W})\partial_{xx} W-(\frac{\partial^2 g}{\partial W^2}\partial_x W)(\partial_y W)-(\frac{\partial g}{\partial W})\partial_{xy} W,\\
\partial_{ty}W &= -(\frac{\partial^2f}{\partial W^2}\partial_y W)(\partial_x W) -(\frac{\partial f}{\partial W})\partial_{xy} W-(\frac{\partial^2 g}{\partial W^2}\partial_y W)(\partial_y W)-(\frac{\partial g}{\partial W})\partial_{yy} W,\\
\partial_{tt}W &= -(\frac{\partial^2 f}{\partial W^2}\partial_t W)(\partial_x W)  -(\frac{\partial f}{\partial W})\partial_{tx} W-(\frac{\partial^2 g}{\partial W^2}\partial_t W)(\partial_y W)-(\frac{\partial g}{\partial W})\partial_{ty} W,
\end{align*}
and so on. These equations can be also used for systems. However, $\frac{\partial f}{\partial W}$ is a matrix and $\frac{\partial^2 f}{\partial W^2}$ is a three-dimensional tensor, etc. The values of all the derivatives at  $\textbf{x}_{G_e}$ are found directly by differentiating the given HWENO reconstruction polynomial with respect to $\textbf{x}$, which is introduced in the following section.

\subsection{HWENO reconstruction in two dimensions}
\label{sechweno2d}

In this subsection, the simple HWENO reconstruction is given on triangular meshes. A direct two dimensional procedure is employed in this paper which is different from the dimension-by-dimension strategy \cite{luo2016,luo2023}. 
To perform the reconstruction,  the cell averages of $W, W_x, W_y$  on the cell $K$ are denoted by
\begin{equation}
\label{aver}
\begin{array}{llll}
|K|\Bar W_{K}&=\int_{K}Wdxdy,\\ 
|K|\Bar V_{K}&=\int_{K}W_xdxdy,\\
|K|\Bar Y_{K}&=\int_{K}W_ydxdy.
\end{array}
\end{equation}

\begin{figure}[hbtp]
\begin{center}
{\includegraphics[width=8cm]{ 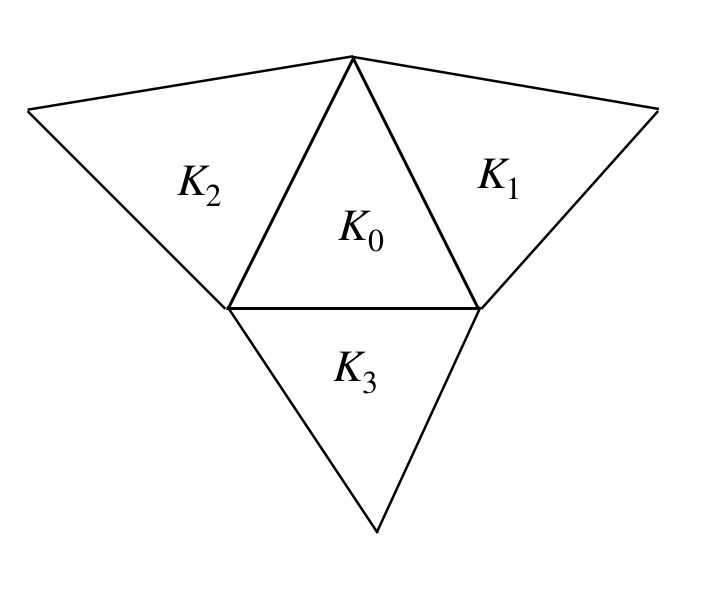}}
\caption{The target cell $K_0$ and its neighboring cells $K_1, K_2, K_3.$ }
\label{figtri}
\end{center}
\end{figure}

Assume $K_0$ is the target cell  and its three neighboring cells are $K_1, K_2, K_3$ as shown in Fig. \ref{figtri} . The reconstruction is as follows.

{\bf Step 1.} Select a big stencil $T_1=\{K_0, K_1, K_2, K_3\}$. Then a quadratic polynomial $p_1(x,y)$ on $T_1$ can be obtained by requiring 
\begin{align*}
\int_{K_{\ell_1}}p_{1}(x,y)dxdy&=|K_{\ell_1}|\Bar W_{K_{\ell_1}},\;\qquad\qquad\ell_1=0,1,2,3,\\
\int_{K_{\ell_2}+K_0}\frac{\partial p_1(x,y)}{\partial x}dxdy&=|K_{\ell_2}|\Bar V_{K_\ell}+|K_0|\Bar V_{K_0},\;\ell_2=1,2,3,\\
\int_{K_{\ell_3}+K_0}\frac{\partial p_1(x,y)}{\partial y}dxdy&=|K_{\ell_3}|\Bar Y_{K_\ell}+|K_0|\Bar Y_{K_0},\;\ell_3=1,2,3,
\end{align*}
which is similar to the work \cite{ji2020}. To obtain a better results, the polynomial $p_1(x,y)$ has the cell averages of $W$ on the target cell $K_0$ and its three neighboring cells while the remaining conditions are matched in a least square sense \cite{hu1999,zhu2018}.

{\bf Step 2.} Choose three small stencils $T_2=\{K_0, K_1\}, T_3=\{K_0, K_2\}$ and $T_4=\{K_0, K_3\}$ and reconstruct three linear polynomials $p_2(x,y), p_3(x,y)$ and $p_4(x,y)$, which satisfy
the following conditions:
\begin{equation*}
	\left\{
	\begin{array}{l}
	\int_{K_0}p_2(x,y)dxdy=|K_0|\Bar W_{K_0},\\
	\int_{K_1}p_2(x,y)dxdy=|K_1|\Bar W_{K_1},\\
	\int_{K_1}\frac{\partial p_2(x,y)}{\partial x}dxdy=|K_1|\Bar V_{K_1},\\
	\int_{K_1}\frac{\partial p_2(x,y)}{\partial y}dxdy=|K_1|\Bar Y_{K_1},
	\end{array}
	\right.
\end{equation*}

\begin{equation*}
	\left\{
	\begin{array}{l}
	\int_{K_0}p_3(x,y)dxdy=|K_0|\Bar W_{K_0},\\
	\int_{K_2}p_3(x,y)dxdy=|K_2|\Bar W_{K_2},\\
	\int_{K_2}\frac{\partial p_3(x,y)}{\partial x}dxdy=|K_2|\Bar V_{K_2},\\
	\int_{K_2}\frac{\partial p_3(x,y)}{\partial y}dxdy=|K_2|\Bar Y_{K_2},
	\end{array}
	\right.
\end{equation*}
and
\begin{equation*}
	\left\{
	\begin{array}{l}
	\int_{K_0}p_4(x,y)dxdy=|K_0|\Bar W_{K_0},\\
	\int_{K_3}p_4(x,y)dxdy=|K_3|\Bar W_{K_3},\\
	\int_{K_3}\frac{\partial p_4(x,y)}{\partial x}dxdy=|K_3|\Bar V_{K_3},\\
	\int_{K_3}\frac{\partial p_4(x,y)}{\partial y}dxdy=|K_3|\Bar Y_{K_3},
	\end{array}
	\right.
\end{equation*}
respectively.

In order to maintain the conservative property, $p_2(x,y), p_3(x,y)$ and $p_4(x,y)$ have the same cell averages as $W$ on the target cell $K_0$ and match the remaining condition in a least square sense \cite{hu1999,zhu2018}.

{\bf Step 3.} Define the linear weights. Employing the similar idea \cite{zhu2016,zhu2017,zhu2018,levy1999,levy2000}, $p_1(x,y)$ is written as
\begin{equation}
\label{pxy}
p_1(x,y)=\gamma_1(\frac{1}{\gamma_1}p_1(x,y)-\sum\limits_{\ell=2}^4\frac{\gamma_\ell}{\gamma_1}p_{\ell}(x,y))+\sum\limits_{\ell=2}^4\gamma_\ell p_{\ell}(x,y)
\end{equation}
Note that \eqref{pxy} holds true for any choice of $\gamma_\ell,\; \ell=1,2,3,4,$ with $\gamma_1 \neq 0$ and $\gamma_\ell, \ell =1,2,3,4$ are the linear weights. Thus the linear weights can be chosen as any positive numbers except the condition that they sum to one. Based on a balance between the essentially non-oscillation in non-smooth regions and accuracy in smooth regions, the linear weights are taken as $\gamma_1=0.94, \gamma_2=\gamma_3=\gamma_4=0.02$ in this paper.

{\bf Step 4.} Compute the smoothness indicators, denoted by $\beta_{\ell}, \ell=1,2,3,4$. To improve the computational efficiency in practice, the same simplified formulation of the smoothness indicators as that in \cite{zhao2025} is adopted, which has the form
\begin{equation}
	\label{smoo1}
	\beta_\ell=\sum\limits_{|s|=1}^k|K_0|^{|s|}(\frac{\partial ^{|s|}p_\ell(x_{K_0},y_{K_0})}{\partial x^{s_1}\partial y^{s_2}})^2, \ell=1,2,3,4,
\end{equation}
where $s=(s_1,s_2), |s|=s_1+s_2$, and for $\ell=1, k=2$; for $\ell=2,3,4, k=1$; $(x_{K_0},y_{K_0})$ is the barycenter of the target cell $K_0$. They measure how smooth the function $p_\ell(x,y)$ is in the target cell $K_0$.

{\bf Step 5.} Compute the nonlinear weights based on the linear weights and smoothness indicators \cite{zhao2025}
\begin{equation}
\label{nonlinearweight}
\omega_l=\frac{\Bar \omega_l}{\Bar\omega_1+\Bar\omega_2+\Bar\omega_3+\Bar\omega_4},\; \Bar \omega_s=\gamma_s(1+\frac{\kappa}{W_{ave}^4\epsilon+W_{ave}^2\beta_s}), l=1,2,3,4; s=1,2,3,4,
\end{equation}
where $W_{ave}=\frac{\sum_{k\in T_1}\Bar{W}_k}{4}+10^{-30},   \kappa=\frac{(|\beta_1-\beta_2|+|\beta_1-\beta_3|+|\beta_1-\beta_4|)^2}{9}$, and $\epsilon$ is a small positive number to avoid the denominator to become zeros, which is set to be $10^{-8}$ in all the computations in this work.

{\bf Step 6.}  Replace the linear weights with the nonlinear weights (\ref{nonlinearweight}), and the final HWENO reconstruction of the conservative values at the point $\text{\bf{x}}_{G_e}$ on the cell $K_0$ is given by
\begin{align}
\label{recon}
W^{int}({\text{\bf{x}}_{G_e}})=\omega_1(\frac{1}{\gamma_1}p_1(\text{\bf{x}}_{G_e})-\sum\limits_{\ell=2}^4\frac{\gamma_\ell}{\gamma_1}p_\ell(\text{\bf{x}}_{G_e}))+\sum\limits_{\ell=2}^4\omega_\ell p_\ell(\text{\bf{x}}_{G_e}).
\end{align}
The derivatives are obtained by differentiating the reconstruction polynomial in  \eqref{recon}. For example, the first derivative with respect to $x$ is given by
 $$\partial_xW^{int}({\text{\bf{x}}_{G_e}})=\omega_1(\frac{1}{\gamma_1}\partial_x p_1(\text{\bf{x}}_{G_e})-\sum\limits_{\ell=2}^4\frac{\gamma_\ell}{\gamma_1}\partial_x p_\ell(\text{\bf{x}}_{G_e}))+\sum\limits_{\ell=2}^4\omega_\ell \partial_xp_\ell(\text{\bf{x}}_{G_e}).$$
The derivatives with respect to the other spatial directions can be obtained in the same manner. The reconstruction of \(W^{\mathrm{ext}}(\mathbf{x}_{G_e})\) and its derivatives is carried out in a similar way.

\subsection{The evaluation of the cell averages of $W_x$ and $W_y$}
As the work \cite{luo2023}, the cell averages of $W_x$ and $W_y$ need to be evaluated at the next time, which is given in the following. Using the divergence theorem in \eqref{aver}, one can obtain 
\begin{equation}
\label{derave2d1}
\begin{array}{llll}
\Bar V_{K}^{n+1}=&\frac{1}{|K|}\int_KW_x(\text{\bf{x}},t_{n+1})dxdy= \frac{1}{|K|}\int_{\partial K}W(\text{\bf{x}},t_{n+1})n_xds,\\
\Bar Y_{K}^{n+1}=&\frac{1}{|K|}\int_KW_y(\text{\bf{x}},t_{n+1})dxdy= \frac{1}{|K|}\int_{\partial K}W(\text{\bf{x}},t_{n+1})n_yds.
\end{array}
\end{equation}
Then the line integrals are approximated by the Gaussian quadrature rule, i.e.,
\begin{equation}
\label{derave2d2}
\begin{array}{llll}
\Bar V_{K}^{n+1}=&\frac{1}{|K|}\int_{\partial K}W(\text{\bf{x}},t_{n+1})n_xds\approx\frac{1}{|K|}\sum\limits_e\sum\limits_{G_e} \hat W(\text{\bf{x}}_{G_e},t_{n+1})n_xw_{G_e}|e|,\\
\Bar Y_{K}^{n+1}=&\frac{1}{|K|}\int_{\partial K}W(\text{\bf{x}},t_{n+1})n_yds\approx\frac{1}{|K|}\sum\limits_e\sum\limits_{G_e} \hat W(\text{\bf{x}}_{G_e},t_{n+1})n_yw_{G_e}|e|,
\end{array}
\end{equation}
where the numerical flux is
\begin{align}
\label{timeflux}
\hat W(\text{\bf{x}}_{G_e},t_{n+1})=\frac{1}{2}(W^{int}(\text{\bf{x}}_{G_e},t_{n+1})+W^{ext}(\text{\bf{x}}_{G_e},t_{n+1})).
\end{align}

In this paper, to get $W^{int}(\text{\bf{x}}_{G_e},t_{n+1})$ and $W^{ext}(\text{\bf{x}}_{G_e},t_{n+1})$, $\tau$ is taken as $\Delta t$ in (\ref{tay}) and (\ref{tayext}) similar to the work \cite{luo2023}, which leads
\begin{align}
\label{taylor2dtime}
W^{int}(\text{\bf{x}}_{G_e},t_n+\Delta t)= W^{int}(\text{\bf{x}}_{G_e},t_n^+)+\sum\limits_{k=1}^{2}[\partial_t^{(k)} W^{int}(\text{\bf{x}}_{G_e},t_n^+)]\frac{(\Delta t)^k}{k!},\\
W^{ext}(\text{\bf{x}}_{G_e},t_n+\Delta t)= W^{ext}(\text{\bf{x}}_{G_e},t_n^+)+\sum\limits_{k=1}^{2}[\partial_t^{(k)} W^{ext}(\text{\bf{x}}_{G_e},t_n^+)]\frac{(\Delta t)^k}{k!},
\end{align}
 where $\partial_t^{(k)} W^{int}(\text{\bf{x}}_{G_e},t_n^+)$ and $\partial_t^{(k)} W^{ext}(\text{\bf{x}}_{G_e},t_n^+)$ have been already computed in (\ref{tay}) and (\ref{tayext}), respectively. Thus the values of the cell averages of $W_x$ and $W_y$ at the next time are evaluated by the Taylor expansion (\ref{derave2d2}). In fact, the Gauss-Lobatto quadrature rule can be employed for the time integration. Thus the values at the next time are calculated and can be directly used in (\ref{timeflux}) for the numerical flux, which can further reduce the computational cost for our method. 

\section{Numerical examples}
\label{secnumexam}
\setcounter{equation}{0}
\setcounter{figure}{0}
\setcounter{table}{0}

\begin{figure}[hbtp]
\begin{center}
{\includegraphics[width=0.5\linewidth]{ 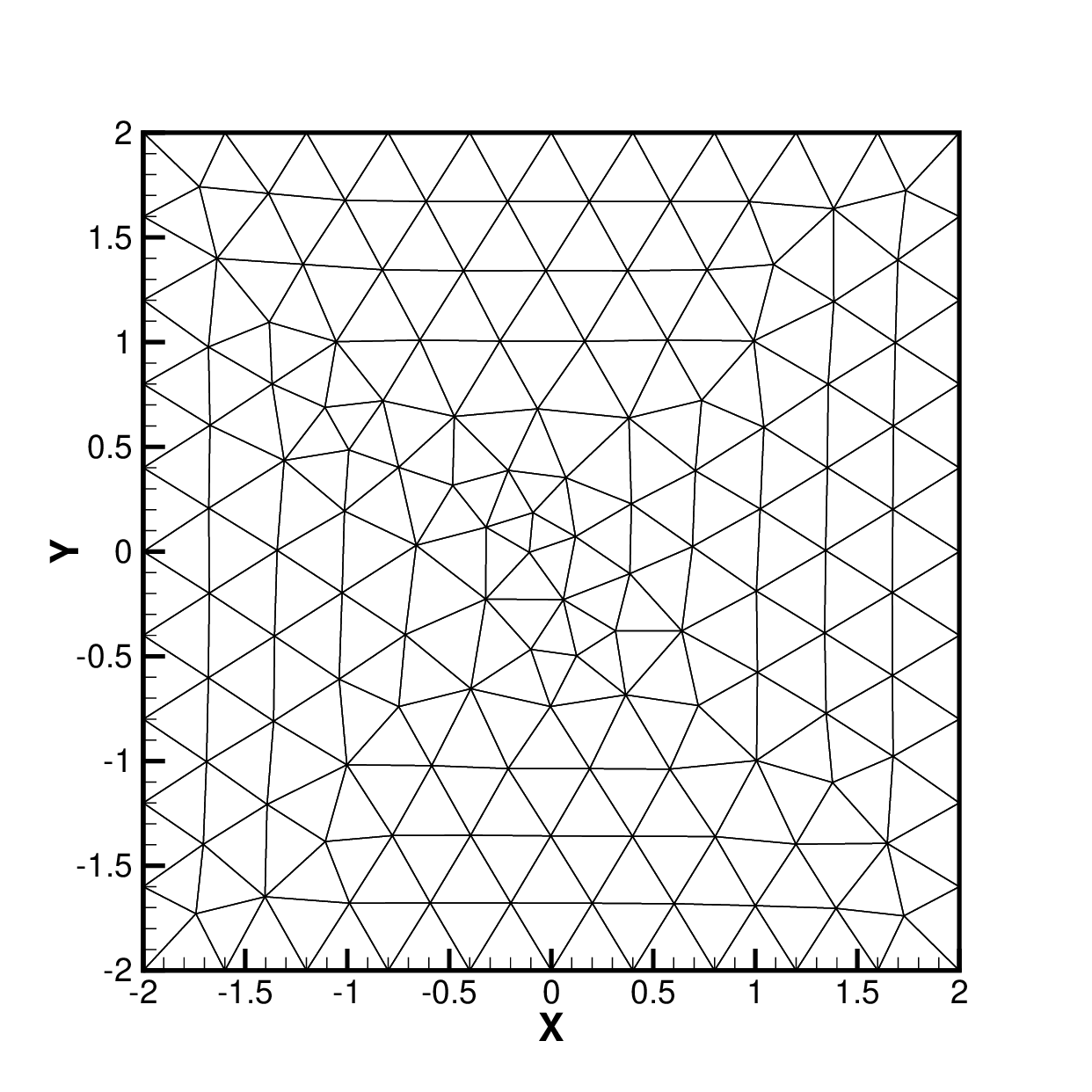}}
\caption{A sample triangular mesh used in two dimensional computation. A  mesh associated with the sample mesh is denoted by $h=\frac{4}{10}$. Triangles: 268. Vertices: 155.}
\label{figbur}
\end{center}
\end{figure}

In this section we present numerical results for a selection of two-dimensional examples to demonstrate the performance of the ADER-HWENO method proposed in the paper. The CFL condition number is taken as 0.5 for all the computations. In addition, we also provide the comparison with the moment-based Hermite WENO schemes \cite{zhao2025}, which is termed as RK-HWENO in the following. Moreover, in order to compare the results with ones in \cite{zhao2025} the same triangular mesh used in their work is employed in the paper; see Fig. \ref{figbur}. A  mesh associated with the initial mesh in Fig. \ref{figbur} is denoted by $h=\frac{4}{10}$. Other meshes are denoted similarly unless otherwise stated. Besides, the $L_1$ and $L_{\infty}$ errors are respectively defined as $\frac{\sum\limits_K|\Bar {W}_K^h-\Bar {W}_K^e|\cdot|K|}{\sum\limits_K|K|}$ and $\max\limits_K|\Bar {W}_K^h-\Bar {W}_K^e|$, where $\Bar {W}_K^h$ represents the numerical cell-average value on the target cell $K$, while $\Bar {W}_K^e$ is the exact cell-average value computed from the exact solution. For most cases, the HWENO reconstruction will perform well. However in some extreme cases, the scale limiter \cite{zhang2011} will be used for modifying the reconstructed polynomials.

\begin{exam}{\em
\label{exam4.1}
We first consider the scalar Burgers' equation in two dimensions
\[
W_t+\left (\frac{W^2}{2}\right )_x+\left (\frac{W^2}{2}\right )_y=0,\quad (x,y)\in (-2,2)\times (-2,2)
\]
subject to the initial condition $W(x,y,0)=0.5+\sin (\frac{\pi(x+y)}{2})$, and a periodic boundary condition in both directions.

We compute the solution up to $T = \frac{0.5}{\pi}$ when the solution is still smooth and the exact solution can be computed using Newton's iteration. The $L_1$ and $L_{\infty}$ errors of the ADER-HWENO method and RK-HWENO method are listed in Table \ref{ex4.1}, which shows the convergence of the third order for both methods. Moreover, from the table one can observe that the errors of the ADER-HWENO method are smaller than those of the RK-HWENO method at the same number of elements.

\begin{table}
\caption{Example~\ref{exam4.1}: Solution error with periodic boundary conditions and $T=\frac{0.5}{\pi}$.}
\renewcommand{\multirowsetup}{\centering}
\begin{center}
\begin{tabular}{|c|c|c|c|c|c|c|c|c|c|c|c|c|}
\hline
method &$N$  &$268$      & $1072$ & $4288$ & $17152$ &$68608$ &274432 \\
\hline
\multirow{4}{2cm}{ADER-HWENO}
 &$L_1$ &6.46e-3  & 6.79e-4 & 7.77e-5 & 1.02e-5 & 1.39e-6 & 1.85e-7 \\
 &  Order & \quad  & 3.25    & 3.13    &  2.93   & 2.87 &2.91 \\
 &$L_{\infty}$ & 3.03e-2 & 6.53e-3 & 8.47e-4  &  1.44e-4 & 2.26e-5&3.34e-6 \\
 &Order    & \quad  & 2.22     & 2.95     &   2.56      & 2.67 &  2.76\\
 \hline
 \multirow{4}{2cm}{RK-HWENO}
 &$L_1$ &1.78e-2  & 3.43e-3 & 2.19e-4 & 1.36e-5 & 1.96e-6 & 2.70e-7 \\
 &  Order & \quad  & 2.37    & 3.97    &  4.02   & 2.78 &2.86 \\
 &$L_{\infty}$ & 7.63e-2 & 2.12e-2 & 3.40e-3  &  1.69e-4 & 3.60e-5&4.94e-6 \\
 &Order    & \quad  & 1.85      & 2.64     &   4.33      & 2.23 &  2.87\\
 \hline
 \end{tabular}
\end{center}
\label{ex4.1}
\end{table}

To see the non-oscillatory properties of the method, the discontinuous solutions computed at $T=\frac{1.5}{\pi}$ with $4288$ triangles along the line $y=x$ are plotted in Fig. \ref{figburg}. From that one can observe the numerical solutions are oscillation-free at discontinuities for all the methods.

\begin{figure}[hbtp]
 \begin{center}
 \mbox{\subfigure[Numerical Solution]
 {\includegraphics[width=0.45\linewidth]{ 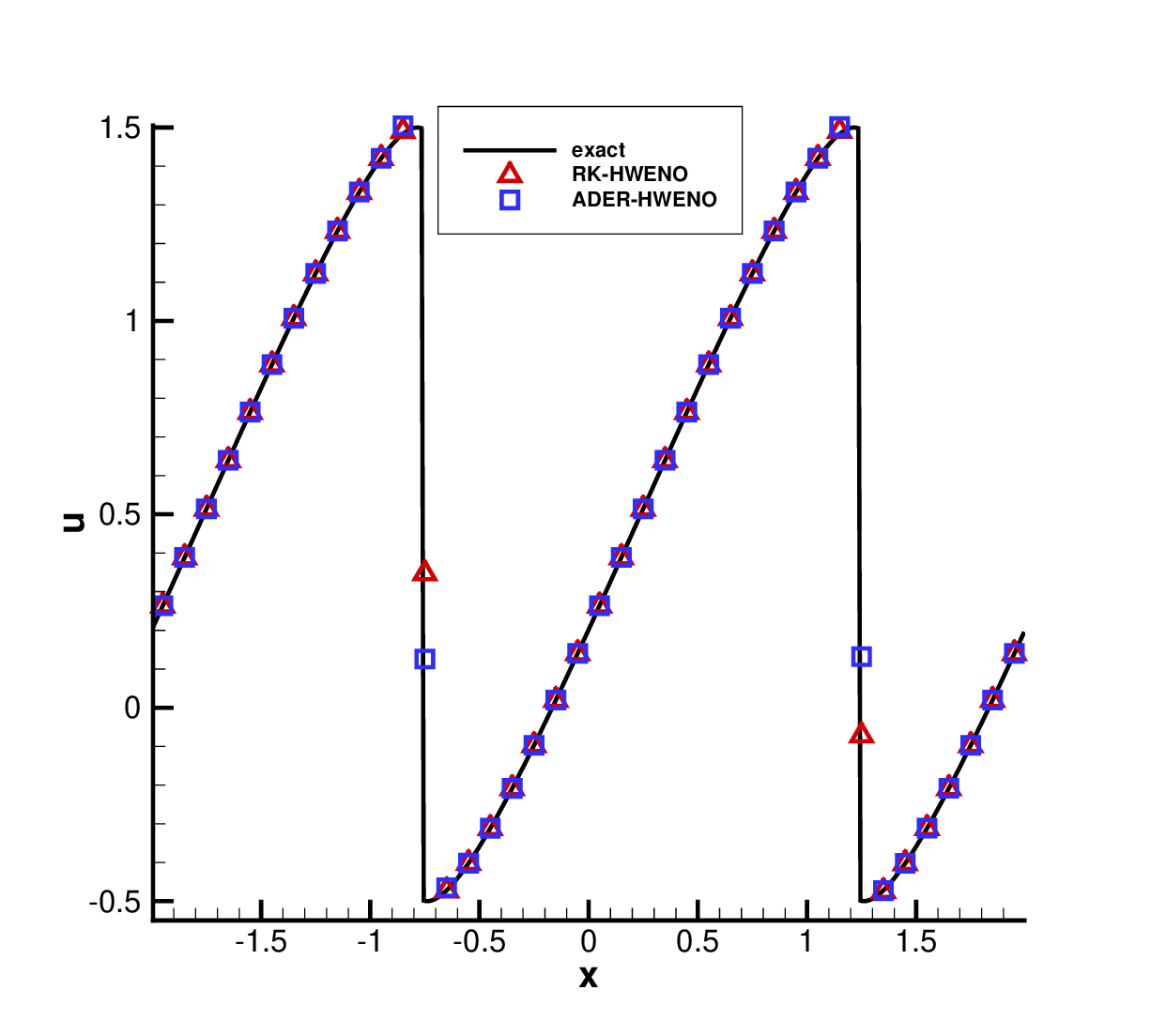}}\quad
 \subfigure[Surface of ADER-HWENO]
 {\includegraphics[width=0.45\linewidth]{ 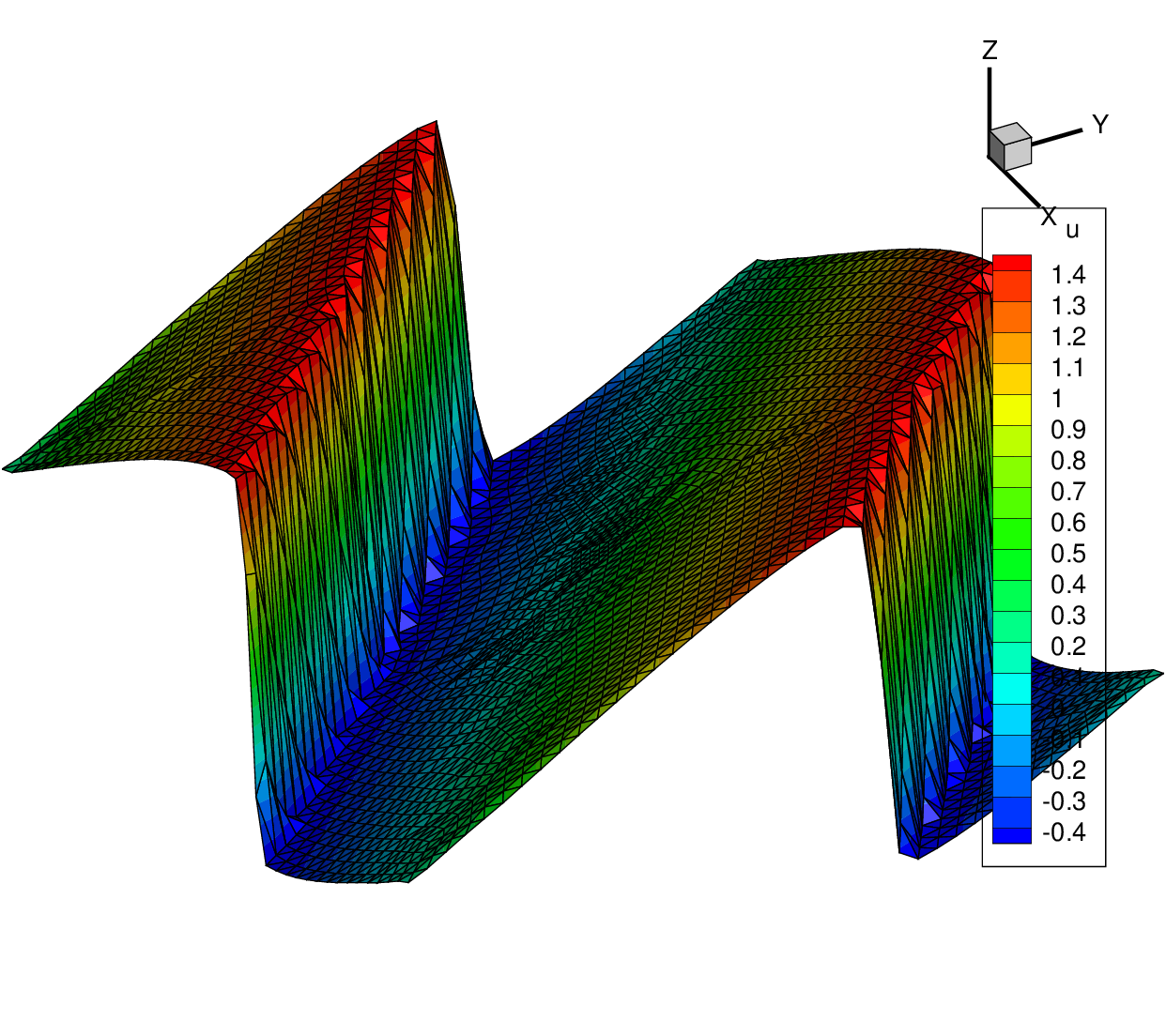}}
   }
   \caption{Example~\ref{exam4.1} Numerical solutions of RK-HWENO and ADER-HWENO method along the line $y=x$. Solid line: the exact solution; Square symbols: ADER-HWENO; Delta symbols: RK-HWENO; Triangles: 4288; Vertices: 2225.}
   \label{figburg}
   \end{center}
   \end{figure}

}\end{exam}

\begin{exam}{\em
\label{exam4.12d}
In order to test the accuracy on the unstructured meshes, we solve the Euler equations
\begin{equation}
\label{2d}
W_t+f(W)_x+g(W)_y\equiv \frac{\partial }{\partial t}
\left(
  \begin{array}{c}
    \rho  \\
     \rho\mu   \\
      \rho\nu   \\
      E \\
   \end{array}
 \right)
 +\frac{\partial }{\partial x}
 \left(
  \begin{array}{c}
    \rho\mu  \\
     \rho\mu^2+P   \\
      \rho\mu\nu   \\
      \mu(E+P) \\
   \end{array}
 \right)
 +\frac{\partial }{\partial y}
 \left(
  \begin{array}{c}
    \rho\nu  \\
     \rho\mu\nu   \\
      \rho\nu^2+P   \\
      \nu(E+P) \\
   \end{array}
 \right)
 =0,
\end{equation}
 where $\rho$ is the density, $\mu$ and $\nu$ are the velocity components
in the $x$- and $y$-direction, respectively, $E$ is the energy density, and $P$ is the pressure. The equation of state is $E=\frac{P}{\gamma-1}+\frac{1}{2}\rho(\mu^2+\nu^2)$ with $\gamma=1.4$.
 The initial condition is given by 
 $$\rho(x,y,0)=1+0.2\hbox{sin}(0.5\pi (x+y)),\; \mu(x,y,0)=1,\; \nu(x,y,0)=1,\; P(x,y,0)=1,$$
  and a periodic boundary condition is applied in both directions.

The computational domain is taken as $(-2,2)\times (-2,2)$ and the stopped time is $T=1$. The results in Table~\ref{ex4.12d} show the convergence of the third order for ADER-HWENO method and RK-HWENO method is achieved for the Euler system in two dimensions. In addition, the $L_\infty$ error of ADER-HWENO method is also smaller than that of RK-HWENO method at the same number of elements.  

To show the efficiency of the ADER-HWENO method, we plot the $L_{\infty}$ error as a function of the CPU time for both ADER-HWENO and RK-HWENO methods in Fig.~\ref{effi}. From the figure, one can see that the ADER-HWENO method is more efficient than the RK-HWENO method. In addition, the CPU time of the ADER-HWENO method is substantially less than that of the RK-HWENO method for the same number of elements.

\begin{table}
	\caption{Example~\ref{exam4.12d}: Solution error with periodic boundary conditions and $T=1$.}
	\begin{center}
		\begin{tabular}{|c|c|c|c|c|c|c|c|c|c|c|c|c|}
			\hline
			method &$N$  &$268$      & $1072$ & $4288$ & $17152$ &$68608$\\
			\hline
			\multirow{4}{1.5cm}{ADER-HWENO}
			&$L^1$       &4.59e-3      & 4.86e-4   & 6.27e-5   & 7.95e-6  & 1.00e-6 \\
			&Order      & \quad        & 3.24     & 2.95     &  2.98        & 2.99 \\
			&$L_{\infty}$ & 1.58e-2      & 2.12e-3  & 1.66e-4   &  1.64e-5 &  2.44e-6\\
			&Order        & \quad         & 2.90     & 3.68      &   3.34      & 2.75   \\
			\hline
			\multirow{4}{1.5cm}{RK-HWENO}
			&$L^1$       &1.67e-2      & 2.14e-3   & 5.81e-5   & 6.26e-6  & 7.93e-7 \\
			&Order      & \quad        & 2.97      & 5.20    &  2.981        & 2.98 \\
			&$L_{\infty}$ & 4.04e-2      & 9.11e-3  & 4.58e-4   &  2.18e-5 &  3.59e-6\\
			&Order        & \quad         & 2.15     & 4.31      &   4.39      & 2.60   \\
			\hline
		\end{tabular}
	\end{center}
	\label{ex4.12d}
\end{table}

\begin{figure}[hbtp]
	\begin{center}
		{\includegraphics[width=0.6\linewidth]{ 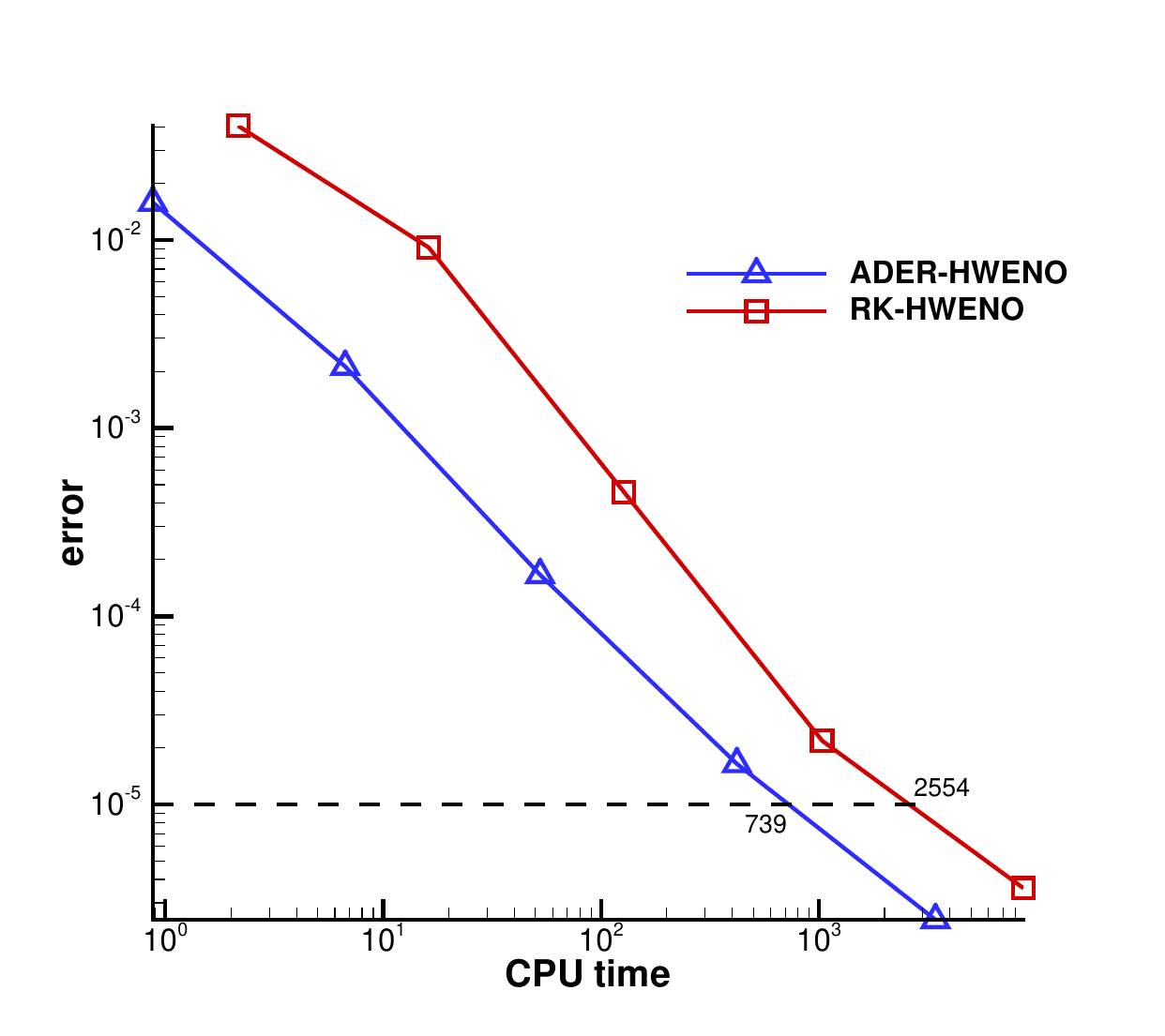}}		
		\caption{$L_{\infty}$ error for both methods is plotted as function of CPU time (in seconds).}
		\label{effi}
	\end{center}
\end{figure}

}\end{exam}

\begin{exam}{\em
\label{isentro}
In this example, the two-dimensional isotropic vortex problem is considered \cite{luo2023,li2021}. The computational domain is taken as $(0,10)\times (0,10)$ and the initial condition is given as follows
\begin{equation*}
\rho= (1-\frac{25(\gamma -1)}{8\gamma \pi^2}e^{1-r^2})^{\frac{1}{\gamma-1}},\\
\mu=1-\frac{5}{2\pi}e^{\frac{1-r^2}{2}}(y-5),\\
\nu=1+\frac{5}{2\pi}e^{\frac{1-r^2}{2}}(x-5),\\
P=\rho^{\gamma},
\end{equation*}
where $r^2=(x-5)^2+(y-5)^2$. The periodic boundary condition is employed in both directions. The exact solution is the vortex along the upper right direction with velocities $(\mu,\nu)=(1,1)$. We compute the numerical solution at the output time $T=1$. The mesh is generated by the GMSH and the sample mesh is given in Fig.~\ref{figisosample}. 

\begin{figure}[hbtp]
	\begin{center}
		{\includegraphics[width=0.6\linewidth]{ 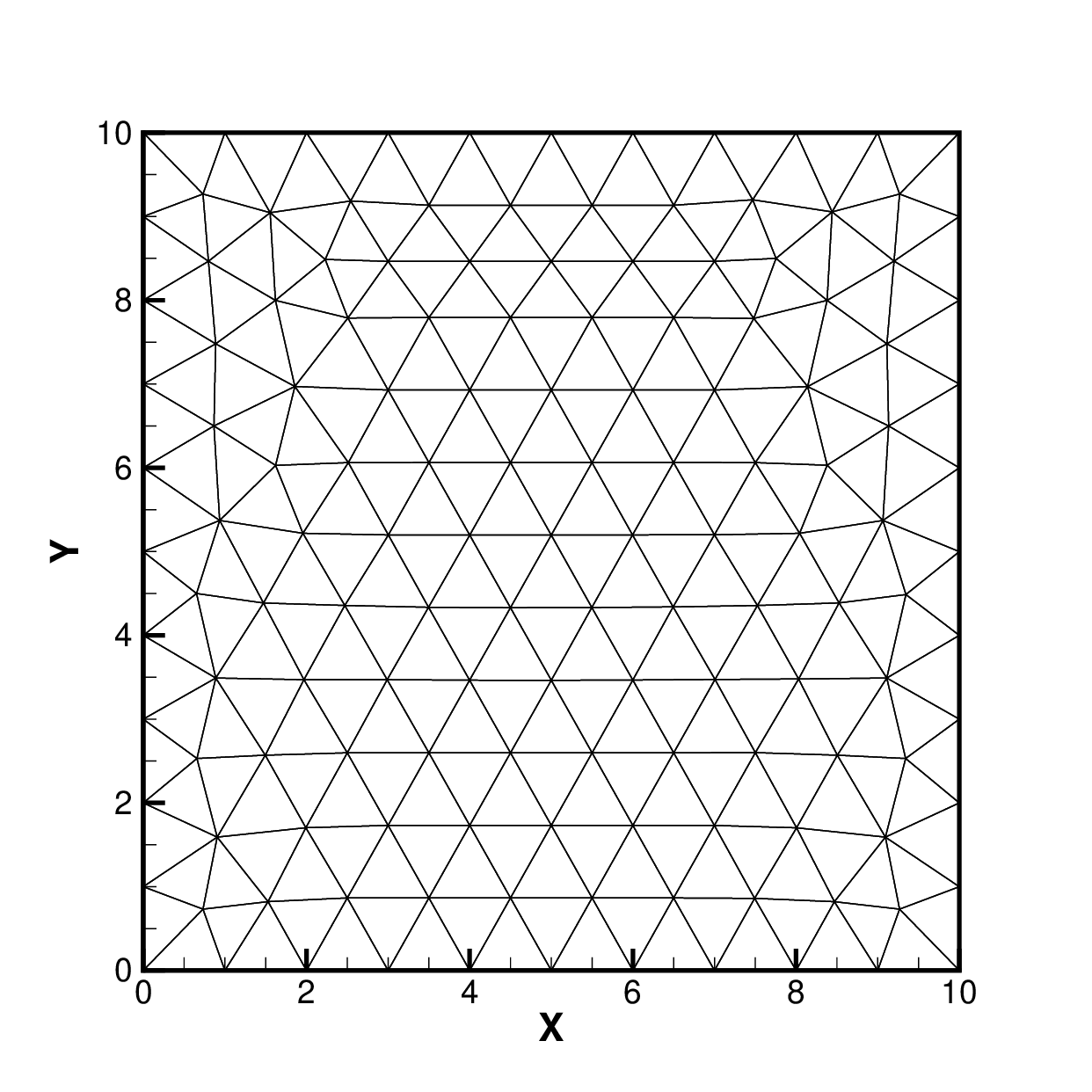}}
		\caption{A sample mesh of the isotropic vortex problem.  The mesh points on the boundary are uniformly distributed with cell length $h=\frac{10}{10}$. Triangles: 244; Vertices:143.}
		\label{figisosample}
	\end{center}
\end{figure}

The errors of the computed density are listed in Table~\ref{isotro} for ADER-HWENO method, in which one can observe that the new method 
can achieve the third order accuracy for this nonlinear smooth problem on the triangular mesh.

\begin{table}
\caption{Example~\ref{isentro}: Solution error with periodic boundary conditions and $T=1$.}
\begin{center}
\begin{tabular}{|c|c|c|c|c|c|c|c|c|c|c|c|c|}
\hline
 $N$  &$244$      & $976$ & $3904$ & $15616$ &$62464$\\
\hline
  $L^1$     &3.45e-3      & 8.16e-4   & 1.33e-4 & 1.76e-5  & 2.22e-6 \\
   Order    & \quad         & 2.08     & 2.62      &  2.92        & 2.98  \\
 $L_{\infty}$ & 5.72e-2   & 1.06e-2  & 1.93e-3  &  2.63e-4 & 3.36e-5\\
 Order       & \quad       & 2.43      & 2.46    &   2.87      & 2.97   \\
 \hline
 \end{tabular}
\end{center}
\label{isotro}
\end{table}

}
\end{exam}

\begin{exam}{\em
\label{examlax}
In this example we consider the Lax problem of the Euler equations (\ref{2d}) with the initial condition
\begin{equation}
(\rho,\mu,\nu,p)=
\left
\{
\begin{array}{ll}
(0.445,0.698,0,3.528), \quad &\text{for}\quad x<0,\\
(0.5,0,0,0.571), \quad &\text{for}\quad x>0.\notag
\end{array}
\right.
\end{equation}
 The inflow/outflow boundary condition is applied in the $x$-direction while the reflective boundary condition is employed in the $y$-direction . The computational domain is $(-0.5,0.5)\times (-0.03,0.03)$ and the integration is stopped at $T=0.16$. 

The problem is simulated by the RK-HWENO and ADER-HWENO on a nonuniform mesh with a triangulation of 100 edges in the $x$-direction and 6 edges in the $y$-direction.
The computed density obtained by these methods along the line $y=0$ is plotted  in Fig. \ref{figlax}, in which one can see that the resolution of the results obtained by ADER-HWENO is comparable with the one obtained by the RK-HWENO method. 

\begin{figure}[hbtp]
 \begin{center}
 \mbox{
 {\includegraphics[width=0.5\linewidth]{ 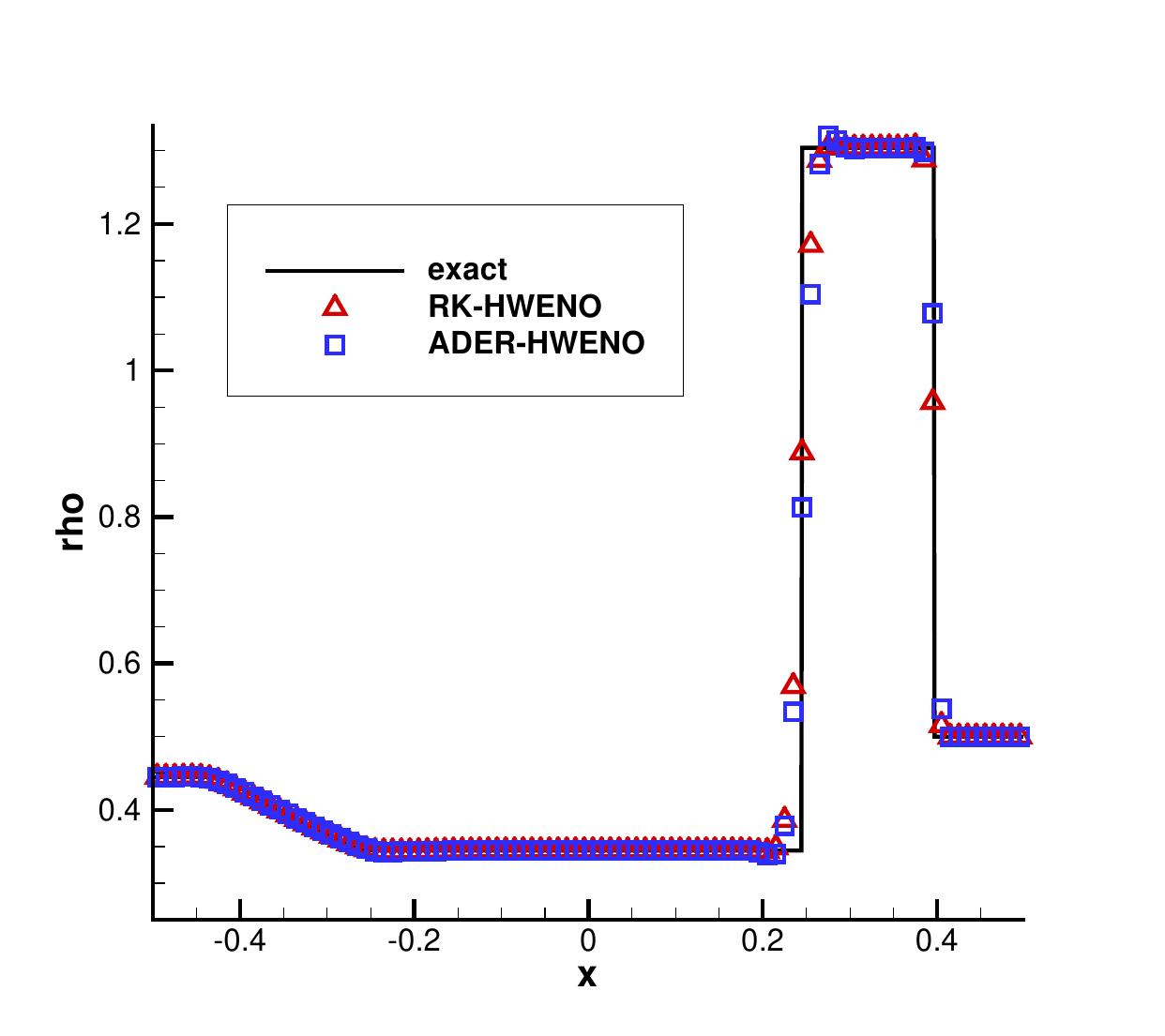}}\quad
 {\includegraphics[width=0.5\linewidth]{ 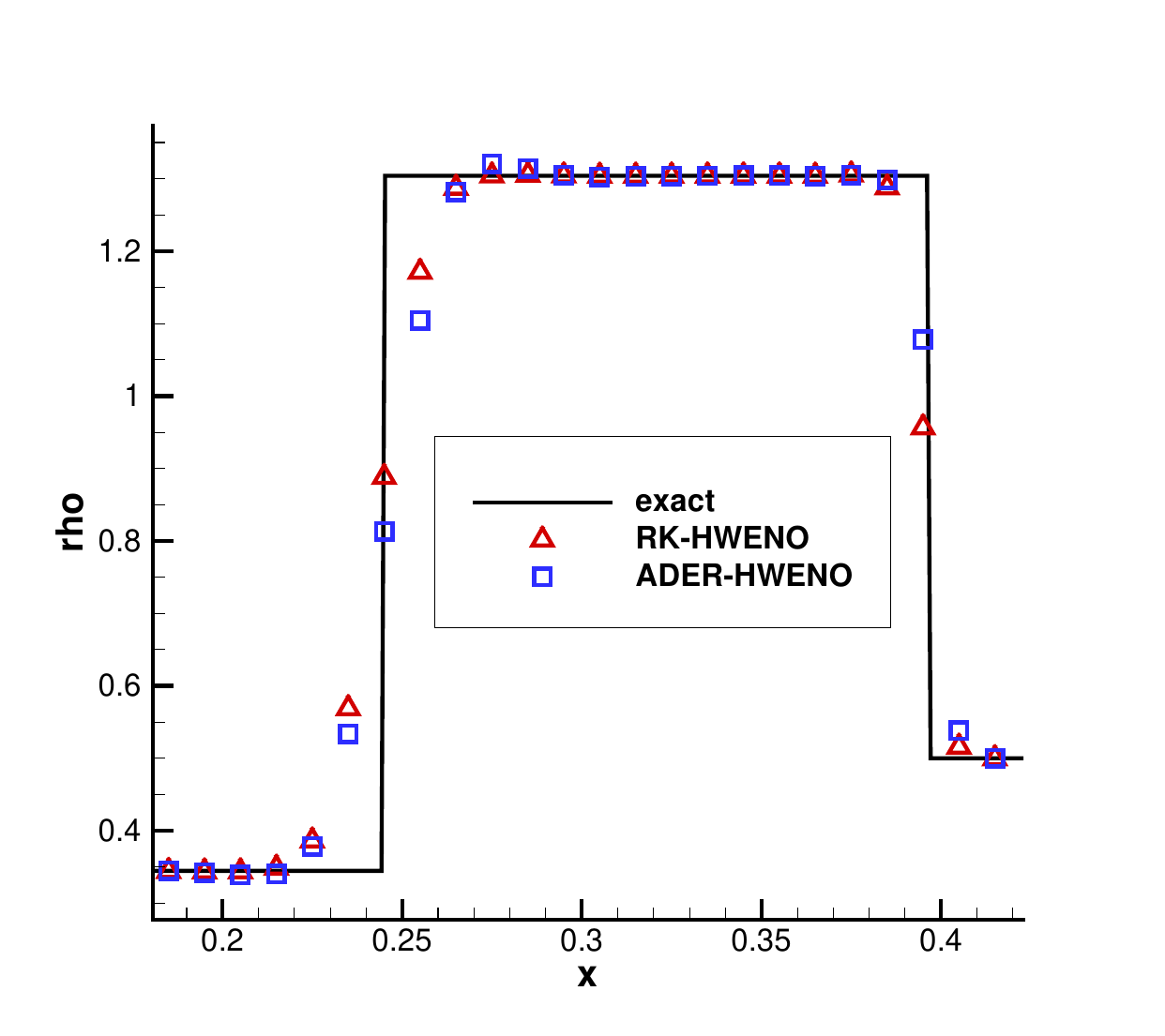}}
}
  \caption{Example~\ref{examlax} The density obtained by RK-HWENO and ADER-HWENO methods. Left: density; Right: zoom of density. Solid line: the exact solution; Square symbols: ADER-HWENO; Delta symbols: RK-HWENO.}
   \label{figlax}
   \end{center}
   \end{figure}

}\end{exam}
\begin{exam}{\em
\label{examshuosher}
The Shu-Osher problem is considered in this example, which contains both shocks and complex smooth region structures. We solve the Euler equations (\ref{2d}) with a moving shock ($\hbox{Mach}=3$) interacting with  a sine wave in density. The initial condition is
\begin{equation}
(\rho,\mu,\nu,p)=
\begin{cases}
(3.857143,2.629369,0,10.333333), \quad &\text{for}\quad x<-4,\\
(1+0.2\hbox{sin}(5x),0,0,1), \quad &\text{for}\quad x>-4.\notag
\end{cases}
\end{equation}
The physical domain is taken as $(-5,5)\times (-0.1,0.1)$ in this computation. The computed density is shown at $T=1.8$ against an ``exact solution'' obtained by a fifth-order finite volume WENO scheme with 10,000 uniform points. 

In this example, a nonuniform mesh with a triangulation of 300 edges in the $x$-direction and 6 edges in the $y$-direction is employed.
The solution obtained by the methods along the line $y=0$ is shown in Fig. \ref{figshuosher}. From the figures, one can observe that the ADER-HWENO method yields higher resolution results than the RK-HWENO methods.

\begin{figure}[hbtp]
\begin{center}
{\includegraphics[width=0.45\linewidth]{ 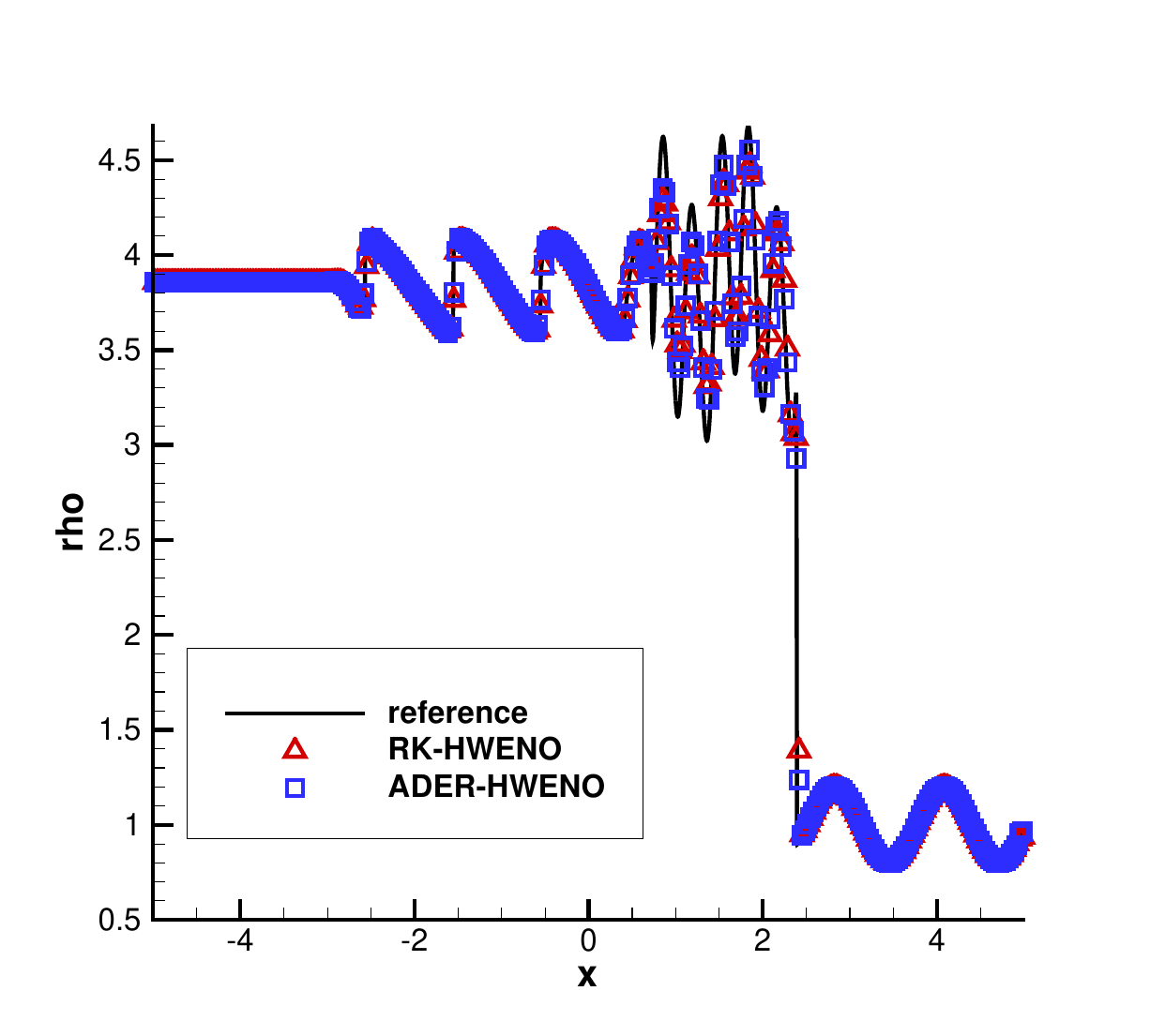}}\quad
{\includegraphics[width=0.45\linewidth]{ 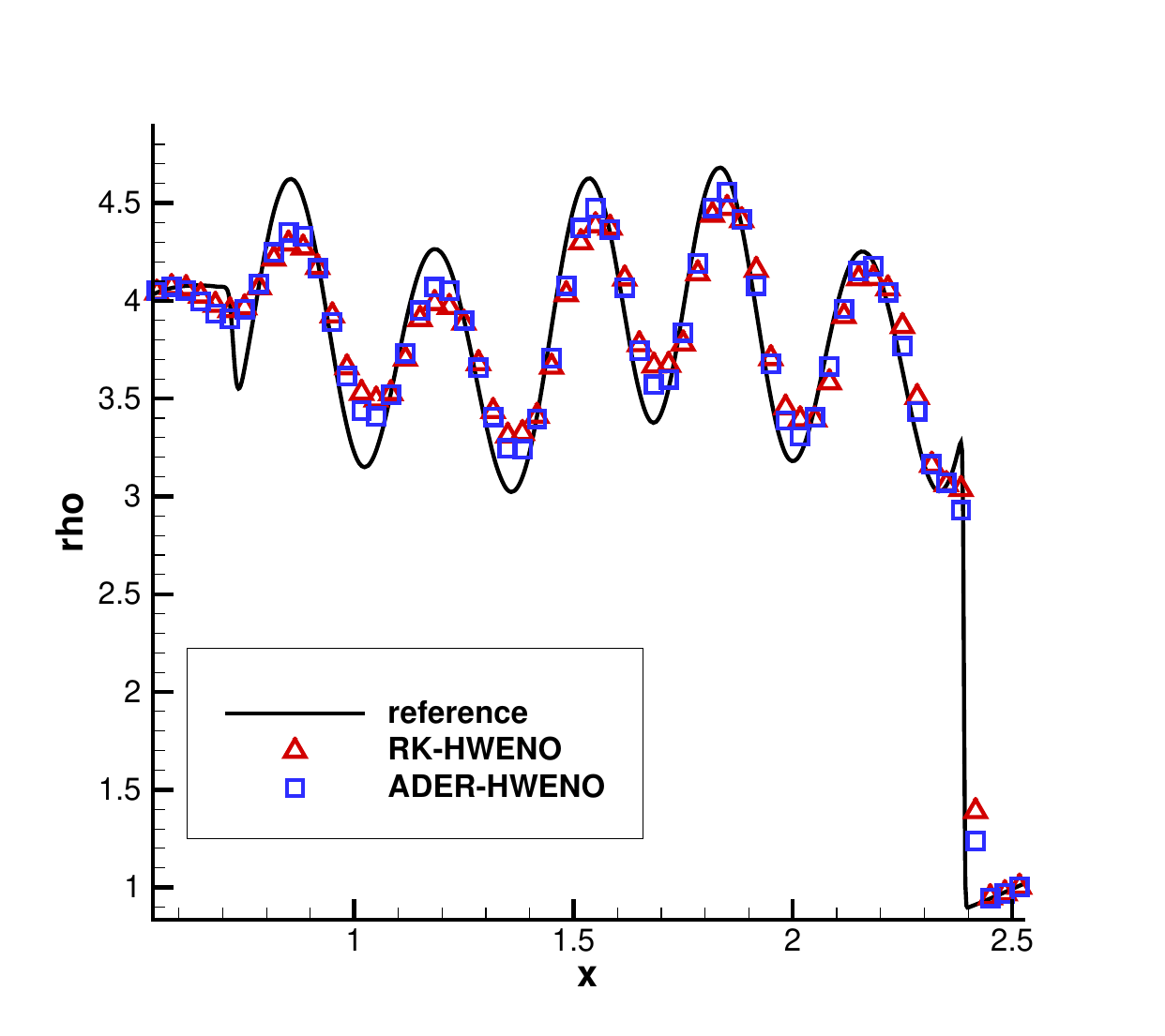}}
\caption{Example~\ref{examshuosher} The density obtained by RK-HWENO and ADER-HWENO methods. Left: density; Right: zoom of density. Solid line: the reference solution; Square symbols: ADER-HWENO; Delta symbols: RK-HWENO. }
\label{figshuosher}
\end{center}
\end{figure}
}
\end{exam}

\begin{exam}{\em
\label{exambw}
We consider the interaction of blast waves of the Euler equations (\ref{2d}), which was first used by Woodward and Colella \cite{woodward1984} as a test problem for various numerical schemes.  The initial condition is given by
\begin{equation}
(\rho,\mu,\nu,p)=
\left
\{
\begin{array}{ll}
(1.0,0,0,1000), \quad& \text{for}\quad 0\leq x<0.1\\
(1.0,0,0,0.01), \quad& \text{for}\quad 0.1\leq x<0.9\\
(1.0,0,0,100),\quad&  \text{for}\quad 0.9\leq x\leq 1 .\notag
\end{array}
\right.
\end{equation}
The physical domain is taken as $(0,1)\times (-0.0075,0.0075)$ and a reflective boundary condition is applied to both ends.

 For the problem, we adopt a nonuniform mesh with a triangulation of 400 edges in the $x$-direction and 6 edges in the $y$-direction. The results at time $T=0.038$ along the line $y=0$ are plotted against an ``exact solution" computed by a fifth-order finite difference WENO scheme \cite{jiang1996} with 81,920 uniform mesh points in Fig. \ref{figbw}. From the figures, one can observe that the ADER-HWENO method provides higher-resolution results than the RK-HWENO method.

\begin{figure}[hbtp]
\begin{center}
{\includegraphics[width=0.45\linewidth]{ 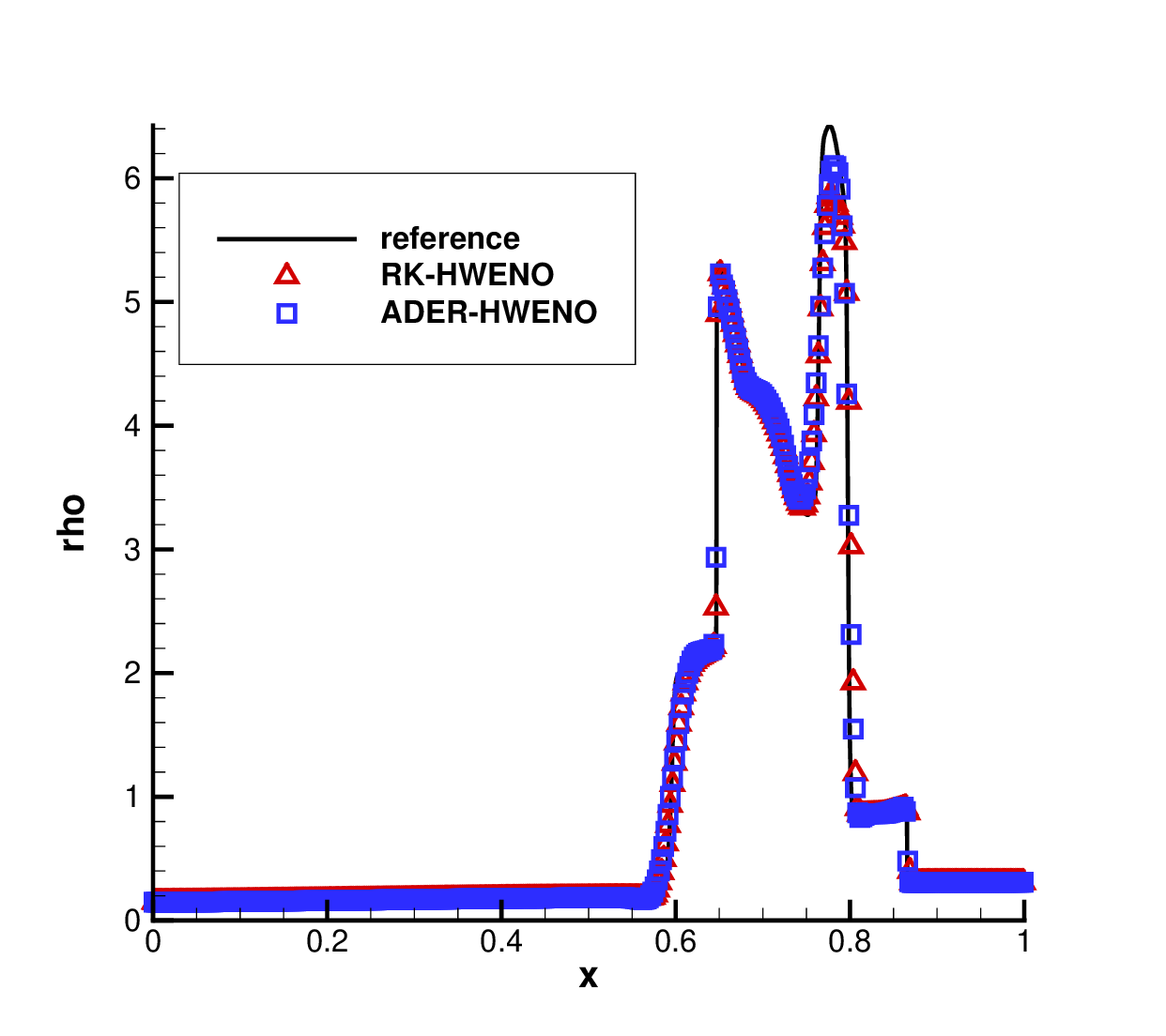}}\quad
{\includegraphics[width=0.45\linewidth]{ 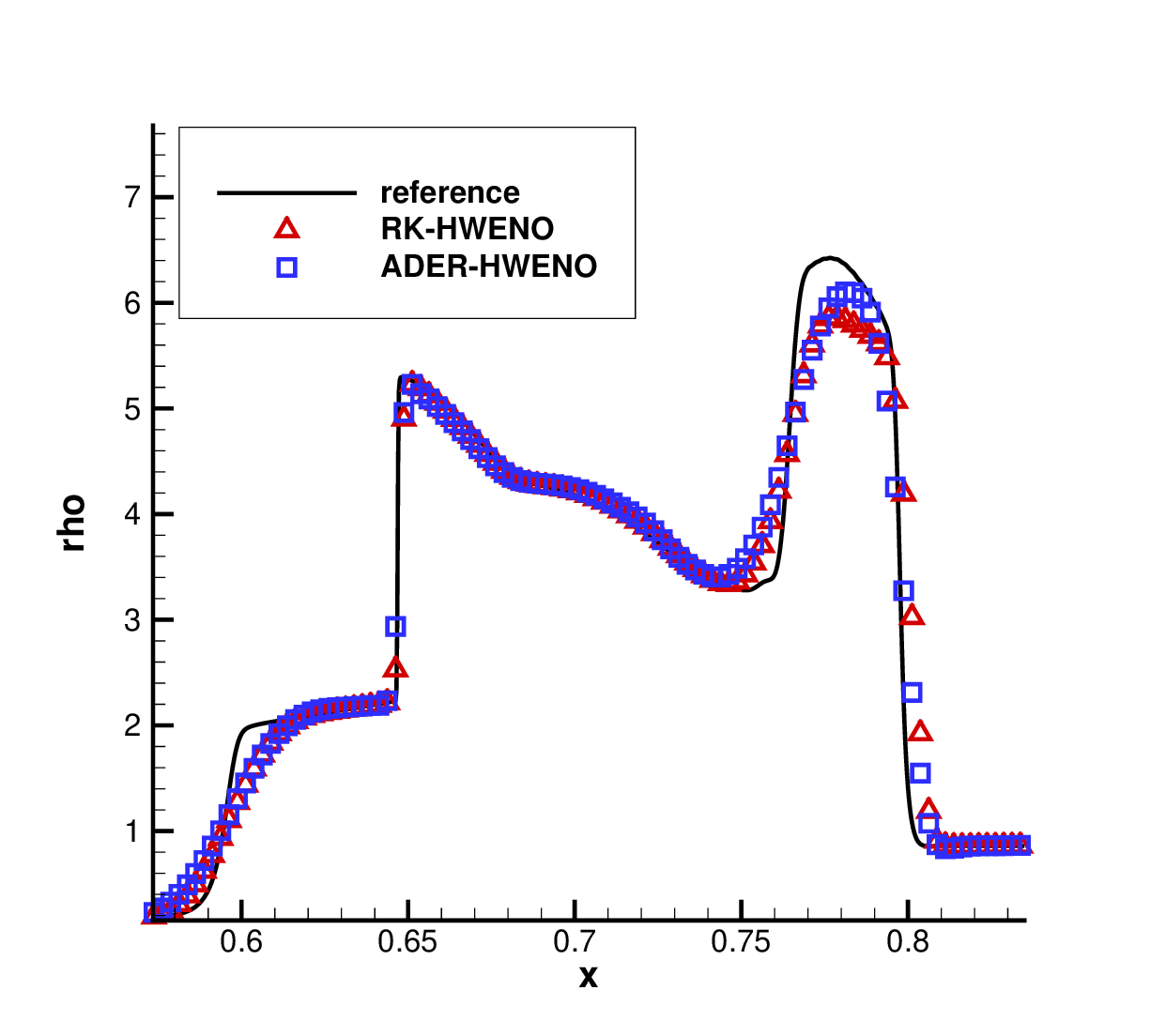}}
\caption{Example~\ref{exambw} The density obtained by RK-HWENO and ADER-HWENO methods. Left: density; Right: zoom of density. Solid line: the reference solution; Square symbols: ADER-HWENO; Delta symbols: RK-HWENO. }
\label{figbw}
\end{center}
\end{figure}
}
\end{exam}

\begin{exam}{\em
\label{doublemach}

This example is the double Mach reflection problem \cite{woodward1984}. We solve the Euler equations (\ref{2d}) in a computational domain of $(0,4)\times (0,1)$. The initial condition is given by 
\begin{equation}
\label{eq1}
W=
\left\{
\begin{array}{ll}
(8,57.1597,-33.0012,563.544)^T, \quad &\text{for} \quad  y\geq h(x,0)\\
(1.4,0,0,2.5)^T, \quad &\text{otherwise}
\end{array}
\right.
\notag
\end{equation}
where $h(x,t)=\sqrt{3}(x-\frac{1}{6})-20t$.
The exact post shock condition is imposed from $0$ to $\frac{1}{6}$ at the bottom while the reflection boundary condition for the rest of the bottom boundary. At the top, the boundary condition is the values that describe the exact motion of the Mach $10$ shock. On the left and right boundaries, the inflow and outflow boundary conditions are used, respectively. The final time is $T=0.2$ using the ADER-HWENO and RK-HWENO method on a mesh with 589,312 triangles and 295,857 vertices. 

The density contours obtained by two methods are shown in Fig. \ref{figdm} on $(0,3)\times (0,1)$, which indicate that the ADER-HWENO and RK-HWENO methods achieve comparable resolution.

\begin{figure}[hbtp]
\begin{center}
{\includegraphics[width=0.8\linewidth]{ 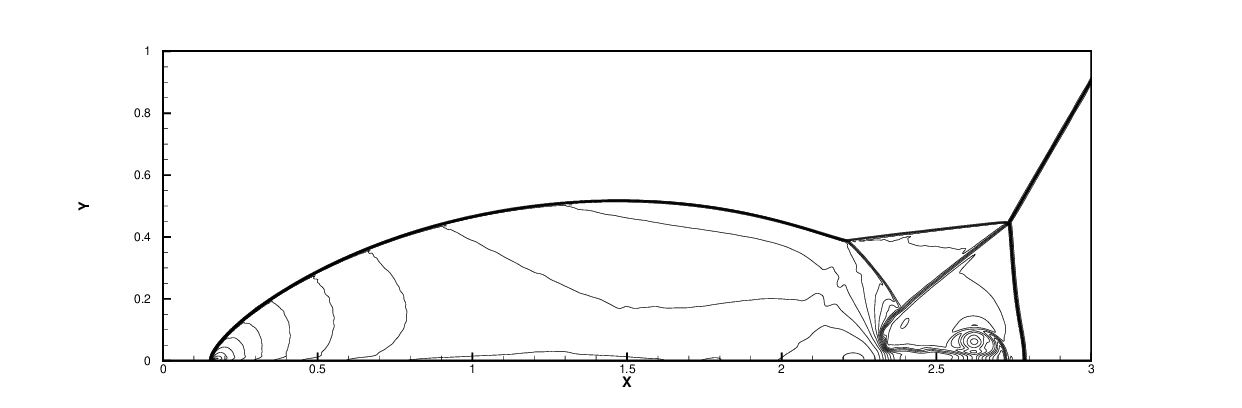}}\\
{\includegraphics[width=0.8\linewidth]{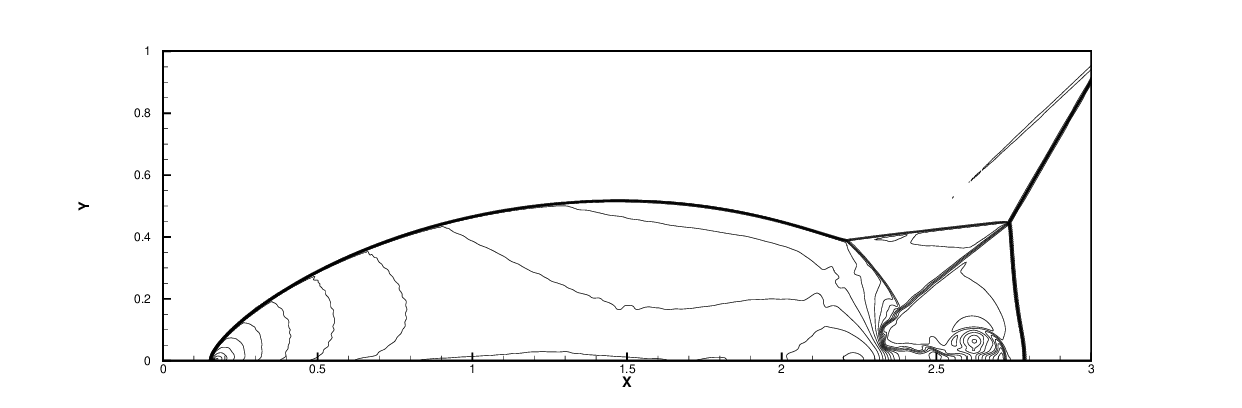}}

\caption{Example~\ref{doublemach} 30 density contours from $1.5$ to 22. Top: ADER-HWENO; Bottom: RK-HWENO. Triangles: 589,312; Vertices: 295,857.}
\label{figdm}
\end{center}
\end{figure}

In order to show the advantages of the unstructured meshes, the other double Mach reflection problem with a more complex domain \cite{woodward1984,zhu2009,ji2020} is considered. And the computational domain is taken as $(-\frac{1}{6},2.5)\times (0,2)$, which contains a wedge $(30^{\circ})$. The initial condition is given by 
\begin{equation}
(\rho,\mu,\nu,p)=
\left\{
\begin{array}{ll}
(8.0,8.25,0.0,116.5)^T, \quad &\text{for} \quad  x<0,\\
(1.4,0,0,1.0)^T, \quad &\text{otherwise}
\end{array}
\right.
\notag
\end{equation}
The exact post shock condition is imposed from  $-\frac{1}{6}$ to $0$ at the bottom. At the top, the boundary condition is the values that describe the exact motion of the Mach $10$ shock. On the left and right boundaries, the inflow and outflow boundary conditions are used, respectively. And the reflection boundary condition for the rest of the boundary. The final time is $T=0.2$. 

In this example, a sample mesh generated by the open software GMSH is shown in Fig. \ref{figdmrsample} denoted by $h=\frac{1}{10}$, which has 852 triangles and 469 vertices. The results using a finer mesh with 209,938 triangles and 837,673 triangles 
obtained by the new method are plotted in Fig. \ref{figdmr}. From the figures, we can observe that the resolution is clearly improved as the mesh is refined.

\begin{figure}[hbtp]
\begin{center}
{\includegraphics[width=0.7\linewidth]{ 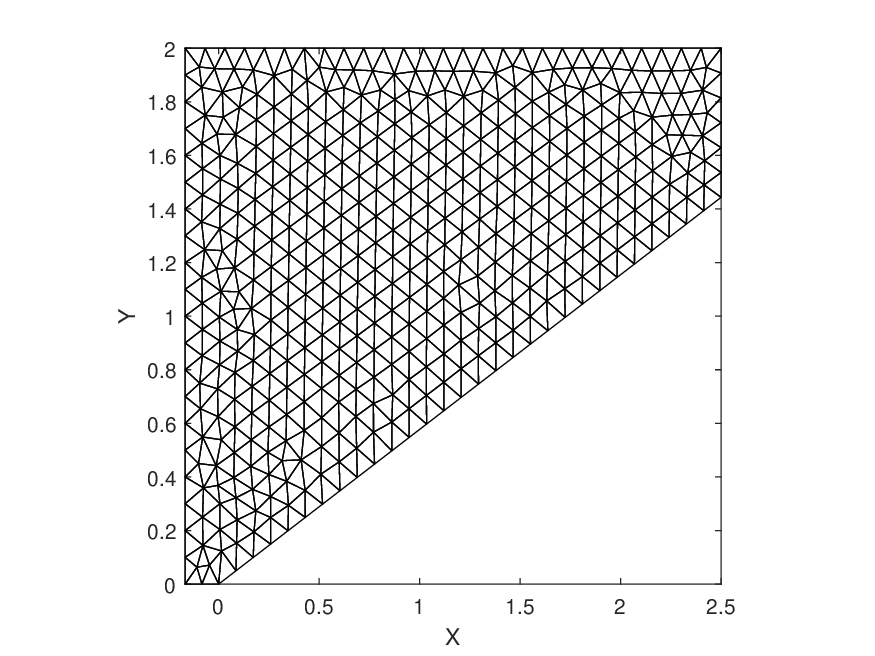}}
\caption{A sample mesh of double mach reflection problem.  The mesh points on the boundary are uniformly distributed with cell length $h=\frac{1}{10}$. Triangles: 852; Vertices: 469.}
\label{figdmrsample}
\end{center}
\end{figure}

\begin{figure}[hbtp]
\begin{center}
\mbox{
{\includegraphics[width=0.5\linewidth]{ 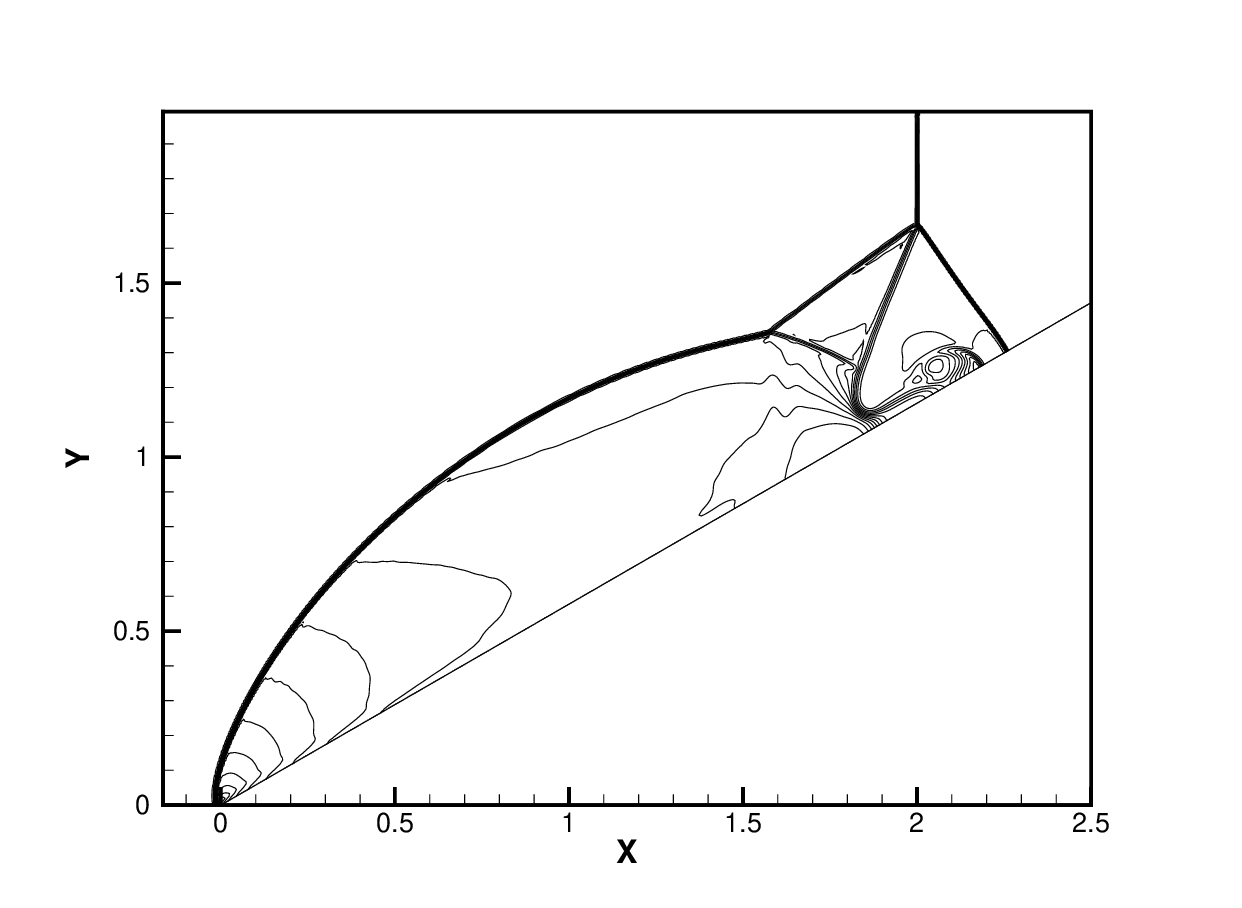}}\quad
{\includegraphics[width=0.5\linewidth]{ 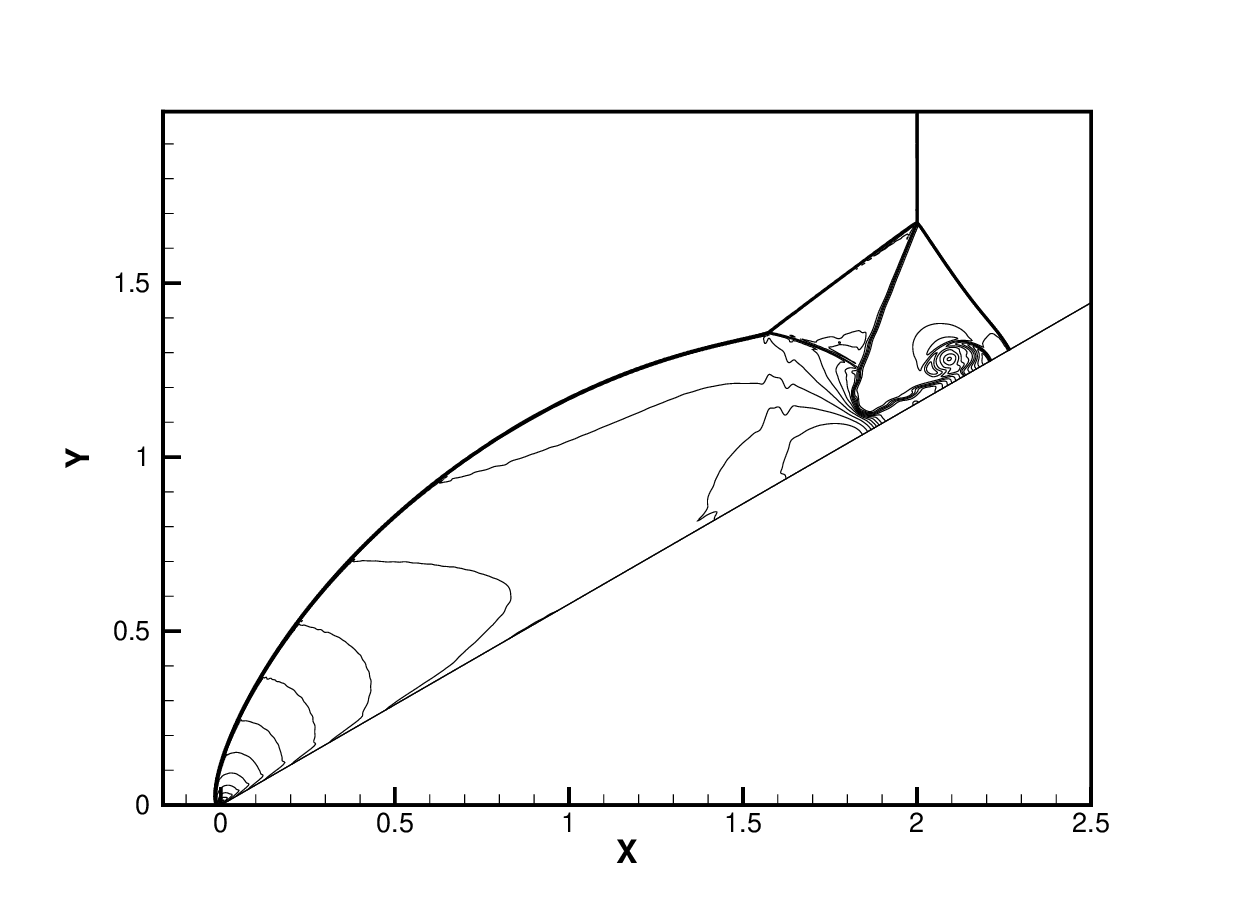}}
}
\mbox{
	{\includegraphics[width=0.5\linewidth]{ 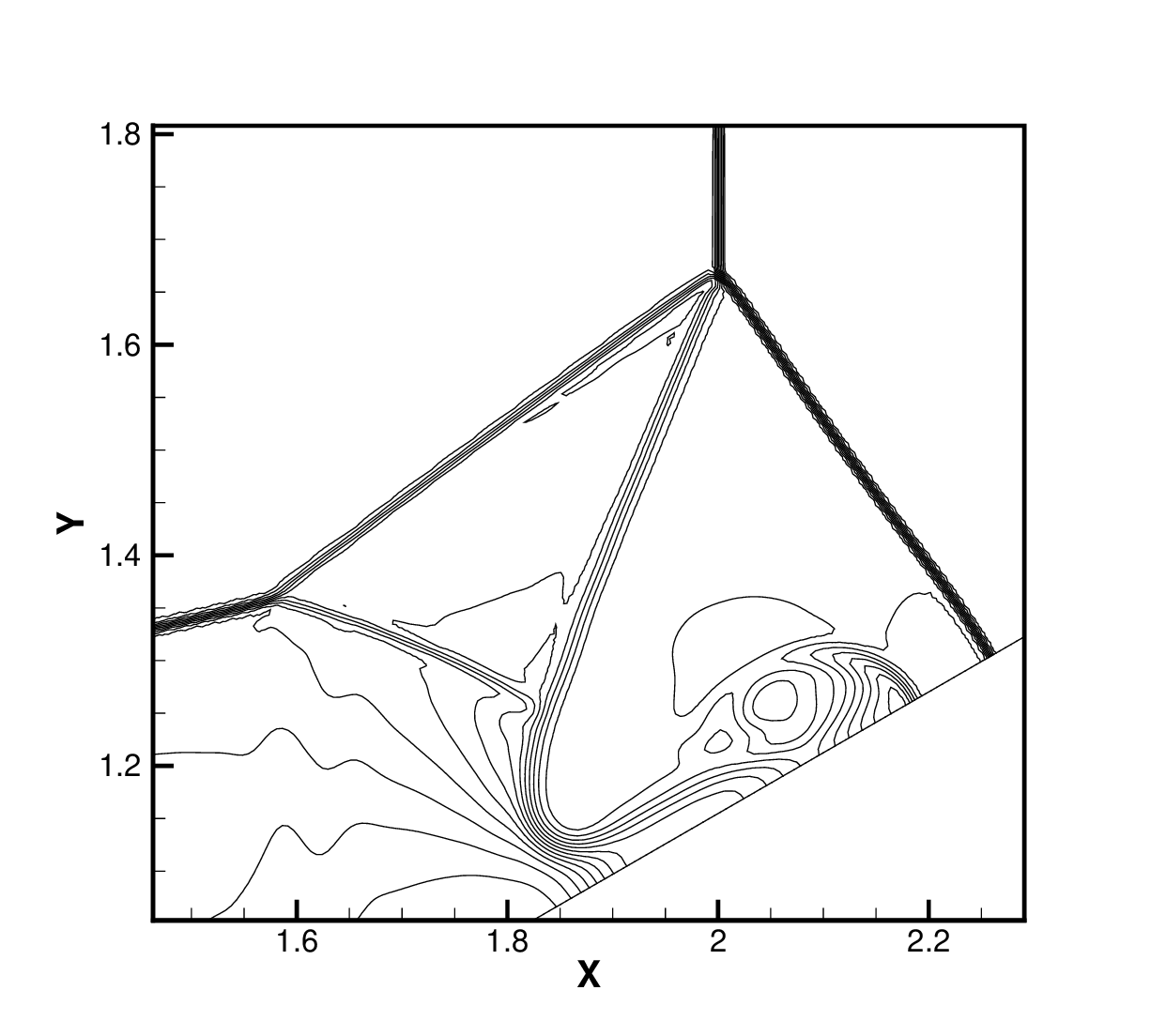}}\quad
	{\includegraphics[width=0.5\linewidth]{ 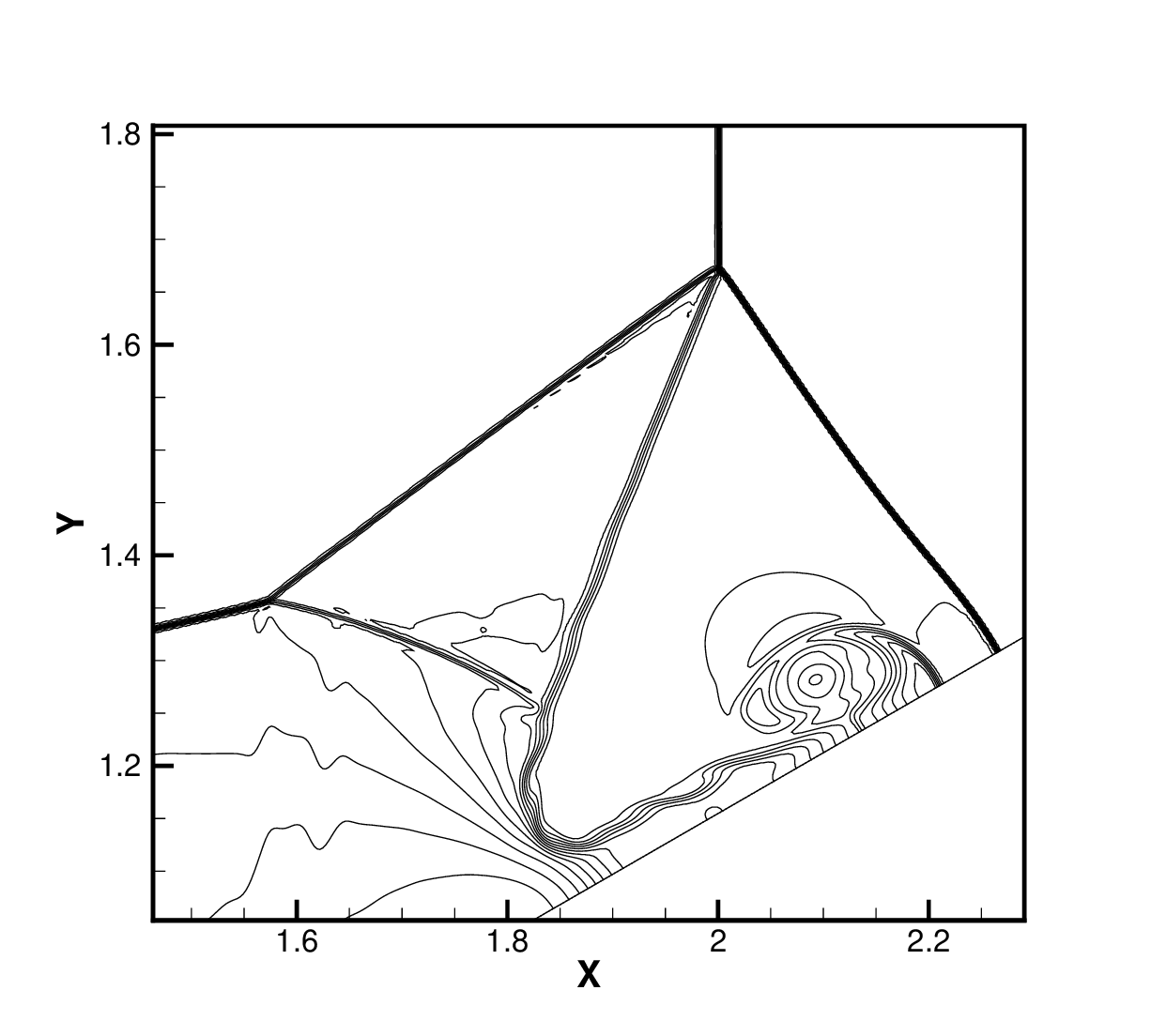}}
}
\caption{Example~\ref{doublemach} 30 density contours from $2$ to 22. Left: $h=\frac{1}{160}$ with 209,938 triangles; Right: $h=\frac{1}{320}$ with 837,673 triangles. Top: the density contours; Bottom: zoom of the complex region.}
\label{figdmr}
\end{center}
\end{figure}

}
\end{exam}

\begin{exam}{\em
\label{step}
 The forward step problem is computed, which is originally from \cite{woodward1984}. We solve the Euler equations (\ref{2d}) in a computational domain of $(0,3)\times (0,1)$. The problem is set up as follows: the wind tunnel is 1 unit wide and 3 units long. The step is 0.2 units high and is located 0.6 units from the left-hand end of the tunnel. The problem is initialized by a right-going Mach 3 flow, namely,
\begin{equation}
(\rho,\mu,\nu,p)=(1.4,3,0,1). \notag
\end{equation}
Reflective boundary conditions are adopted along the wall of the tunnel and inflow and outflow boundary conditions are applied at the entrance and exit, respectively. For this problem a singularity is located at the corner of the step. However, we do not modify the scheme or refine the mesh near the corner to test the new method, which is different from the work \cite{ji2020}. The final time is $T=4$. 

 A sample mesh for the example is shown in Fig. \ref{figstepsample}. The density computed by the ADER-HWENO and RK-HWENO methods on a refiner mesh with 164,608 triangles and 82,945 vertices are plotted in Fig. \ref{figstep}. From the figures, we can see that the resolution of the new method is comparable with that of RK-HWENO method. 

\begin{figure}[hbtp]
\begin{center}
{\includegraphics[width=0.8\linewidth ]{ 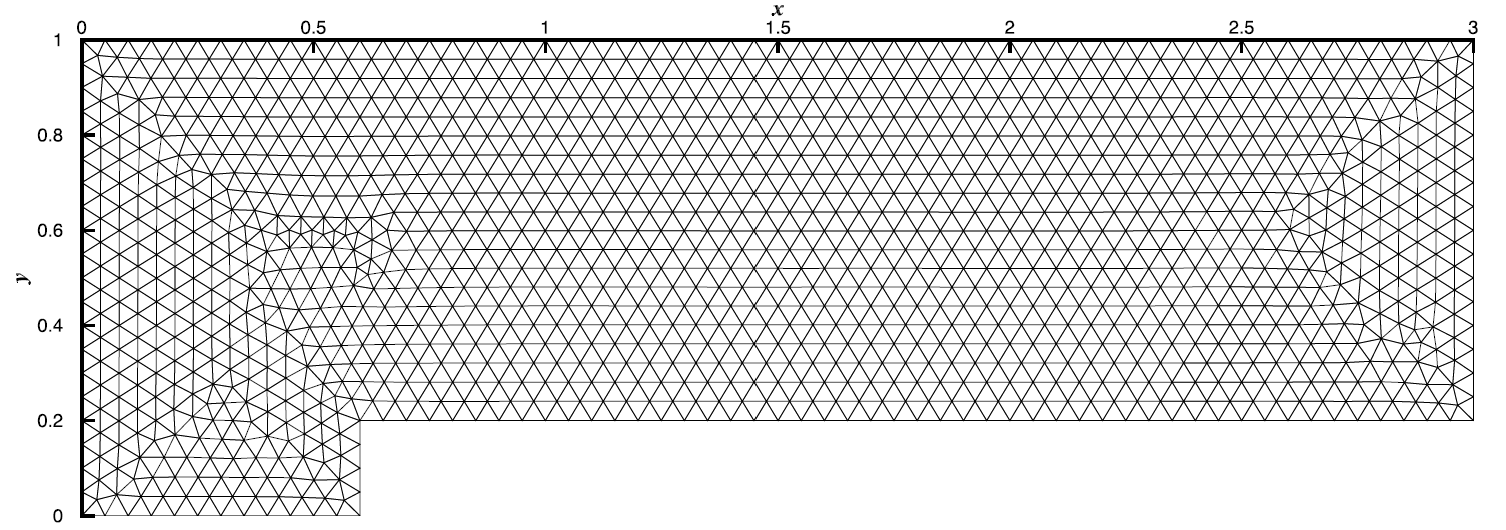}}
\caption{A sample mesh of forward step problem.  The mesh points on the boundary are uniformly distributed with cell length $h=\frac{1}{20}$. Triangles: 2572; Vertices: 1367.}
\label{figstepsample}
\end{center}
\end{figure}

\begin{figure}[hbtp]
\begin{center}
{\includegraphics[width=0.8\linewidth]{ 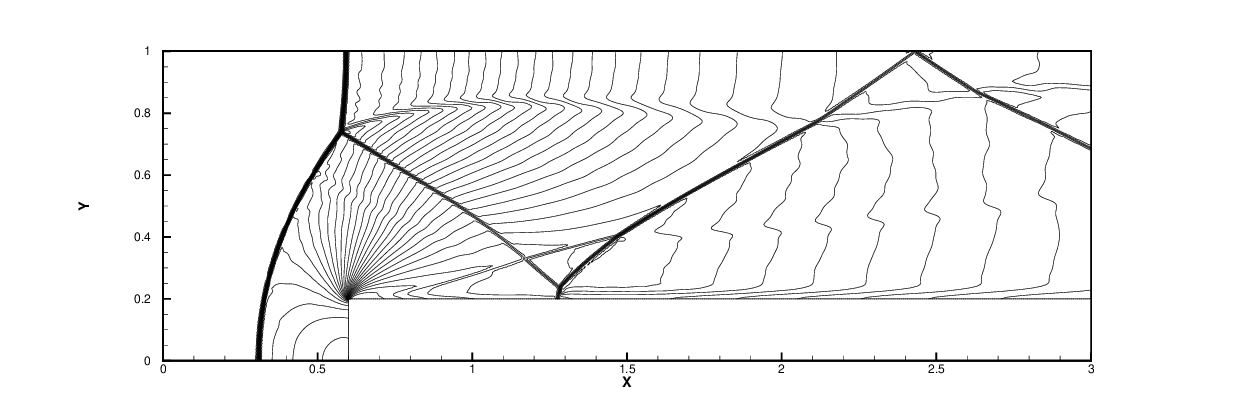}}\quad
{\includegraphics[width=0.8\linewidth]{ 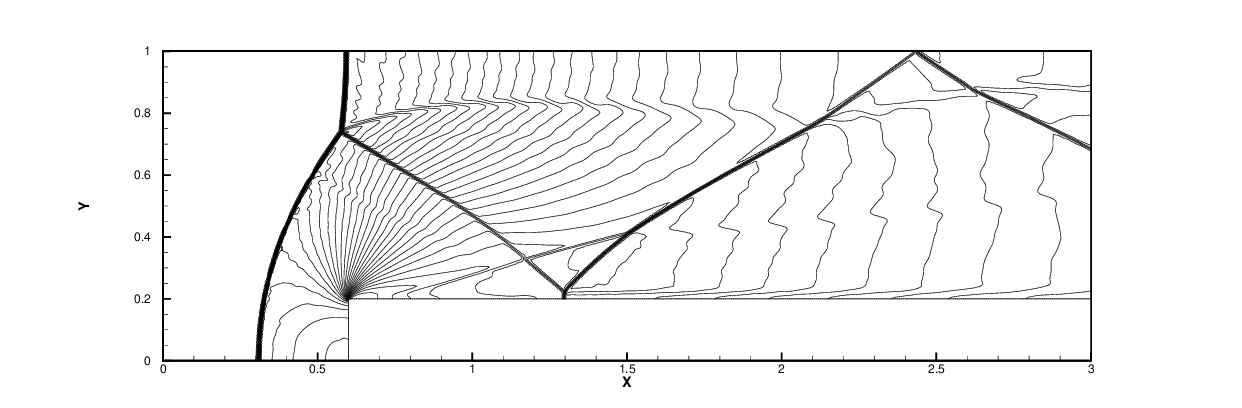}}
\caption{Example~\ref{step} 30 density contours from $0.32$ to 6.15. Top: ADER-HWENO; Bottom: RK-HWENO.}
\label{figstep}
\end{center}
\end{figure}
}
\end{exam}

\section{Conclusions}
\label{secconclusion}
\setcounter{equation}{0}
\setcounter{figure}{0}
\setcounter{table}{0}

We have presented a compact and high order HWENO scheme using the ADER time discretization for hyperbolic conservation laws on the triangular mesh. The Lax-Wendroff procedure is adopted to convert time derivatives to spatial derivatives. Benefiting from this, the cell averages of the derivatives are obtained by the Gaussian points along the cell interfaces through the Green-Gauss theorem instead of by the evolution solution directly. Compare with the RK-HWENO scheme in \cite{zhao2025} on the triangular meshes, the method newly developed in this paper has several advantages. Firstly, RK-HWENO method must solve the additional equations for reconstructions and time advancing, which is avoided for the new method. Secondly, the HWENO reconstruction of the ADER-HWENO method is performed once every time step and is different from the RK-HWENO method, in which the reconstruction need to be performed several times for a time step. Thus the ADER-HWENO method is more efficient than the RK-HWENO method, which is also observed in the example of accuracy test. For the same order of accuracy, the new method employs a more compact reconstruction stencil than the ADER-WENO method \cite{dumbser20071,dumbser20072}, since the HWENO formulation evolves both the solution and its derivative moments. Numerical results indicate the high order for smooth solutions both in space and time and keep non-oscillatory at discontinuities. We recall that the work is in the Euler framework for hyperbolic conservation laws. Extending the developed method to the moving mesh \cite{luo2019,luo2022} and multi-component flows \cite{luo2021} has been underway.
 
\section*{Acknowledgements}
{
This work is supported partly by National Key R\&D Program of China (Grant Number 2022YFA1004500), National Natural Science Foundation of China (Grant Nos. 12571416, 12271052), the Foundation of National Key Laboratory of Computational Physics (Grant Number 6142A05240306), National Natural Science Foundation of Xiamen, China (Grant Number 3502Z202472004), and Fujian Provincial Natural Science Foundation of China (Grant Number 2026J009009).
}


\end{document}